\documentclass{article}
\usepackage{arxiv}
\usepackage[utf8]{inputenc} 
\usepackage[T1]{fontenc}    
\usepackage{hyperref}       
\usepackage{url}            
\usepackage{booktabs}       
\usepackage{amsfonts}       
\usepackage{nicefrac}       
\usepackage{microtype}      
\usepackage{lipsum}		
\usepackage{float}
\usepackage{graphicx}
\usepackage{natbib}
\usepackage{doi}
\usepackage{float}
\usepackage{amsmath}
\usepackage{todonotes}

\DeclareMathOperator*{\argmin}{argmin}

\usepackage{soul}

\title{Data-driven modelling of brain activity using neural networks, Diffusion Maps, and the Koopman operator}
\author{ 
\hspace{1mm}Ioannis K.~Gallos\\
	School of Applied Mathematics and Physics\\
	National Technical University of Athens\\
	Athens, Greece \\
\texttt{yiannis.gallos@gmail.com}
	\And
\hspace{1mm}Daniel Lehmberg\\
	School of Computation, Information and Technology\\
	Technical University of Munich\\
	Munich, Germany \\
	\texttt{d.lehmberg@tum.de} \\
	\And
\hspace{1mm}Felix Dietrich\\
	School of Computation, Information and Technology\\
	Technical University of Munich\\
	Munich, Germany \\
	\texttt{felix.dietrich@tum.de} \\
	\And
\hspace{1mm}Constantinos ~Siettos\thanks{Corresponding author: constantinos.siettos@unina.it}\\
	Department of Mathematics and Applications\\
	University of Naples Federico II\\
	Naples, Italy \\
\texttt{constantinos.siettos@unina.it} \\
}

\renewcommand{\shorttitle}{\textit{arXiv} Template}
\hypersetup{
pdftitle={Data-driven modelling of brain activity using machine learning and Koopman operator},
pdfsubject={q-bio.NC, q-bio.QM},
pdfauthor={Ioannis K.~Gallos, Daniel Lehmberg, Felix Dietrich, Constantinos I.~Siettos},
pdfkeywords={},
}
\begin{document}
\maketitle

\begin{abstract}
We propose a machine-learning approach to model long-term out-of-sample dynamics of brain activity from task-dependent fMRI data. Our approach is a three stage one. First, we exploit Diffusion maps (DMs) to discover a set of variables that parametrize the low-dimensional manifold on which the emergent  high-dimensional fMRI time series evolve. Then, we construct reduced-order-models (ROMs) on the embedded manifold via two techniques: Feedforward Neural Networks (FNNs)  and the Koopman operator. Finally, for predicting the out-of-sample long-term dynamics of brain activity in the ambient fMRI space, we solve the pre-image problem coupling DMs with Geometric Harmonics (GH) when using FNNs and the Koopman modes per se. For our illustrations, we have assessed the performance of the two proposed schemes using a benchmark fMRI dataset with recordings during a visuo-motor task. The results suggest that just a few (for the particular task, five) non-linear coordinates of the high-dimensional fMRI time series provide a good basis for modelling and out-of-sample prediction of the brain activity. Furthermore, we show that the proposed approaches outperform the one-step ahead predictions of the naive random walk model, which, in contrast to our scheme, relies on the knowledge of the signals in the previous time step. Importantly, we show that the proposed Koopman operator approach provides, for any practical purposes, equivalent results to the FNN-GH approach, thus bypassing the need to train a non-linear map and to use GH to extrapolate predictions in the ambient fMRI space; one can use instead the low-frequency truncation of the DMs function space of $L^2$-integrable functions, to predict the entire list of coordinate functions in the fMRI space and to solve the pre-image problem.
\end{abstract}

\section{Introduction}
Understanding and modelling the emergent brain dynamics has been a primary challenge in contemporary neuroscience \citep{bullmore2009complex,ermentrout2010mathematical,jorgenson2015brain,siettos2016multiscale,breakspear2017dynamic,sip2023characterization}. Towards this goal, different approaches at a variety of scales have been introduced, ranging from individual neuron dynamics to the behaviour of different regions of neurons across the brain \citep{izhikevich2007dynamical,deco2008dynamic,bullmore2012economy,papo2014complex,siettos2016multiscale,breakspear2017dynamic,sip2023characterization}. Some of them are physics-biologically informed models, thus derived directly from first principles (such as the Hodgkin–Huxley \citep{hodgkin1952propagation,nelson1998hodgkin,mccormick2007hodgkin,spiliotis2022deep}, the FitzHugh-Nagumo \citep{fitzhugh1955mathematical,nagumo1962active}, phase-oscillators \citep{kopell1986symmetry,laing2017phase,skardal2020higher}) and neural mass models \citep{segneri2020theta,taher2020exact}, the Baloon model \citep{buxton1998dynamics} and Dynamic Causal Modelling (DCM) \citep{penny2004modelling}) (for a review of the various dynamic models see also \citep{izhikevich2007dynamical}). On the other hand, there is the data-driven approach, exploiting a wide range of methods (see also \citep{heitmann2018brain}) extending from independent component analysis \citep{beckmann2004probabilistic}, to Granger causality-based models \citep{seth2010matlab,seth2015granger,friston2013analysing,protopapa2014granger,kugiumtzis2015direct,protopapa2016children,kugiumtzis2017dynamics,almpanis2020construction}, phase-synchronization \citep{mormann2000mean,rudrauf2006frequency,jirsa2013cross,zakharova2014chimera,mylonas2016modular,scholl2022partial} and network-based ones \citep{kopell2004chemical,spiliotis2022complex,petkoski2019transmission}. For an extended  review of the various multiscale approaches, see \citep{siettos2016multiscale}.\par
Most of the multiscale data-driven modelling approaches are ``blended'', i.e., they rely on EEG, MEG and fMRI data to calibrate biological-inspired models or to build surrogate models for the emergent brain activity. Within, this framework, machine learning (ML) has also been exploited to build surrogate dynamical models from neuroimaging data \citep{grossberg1992neural,niv2009reinforcement,richiardi2013machine,suk2016state,gholami2018modelling,sun2022deep}. However, when it comes to fMRI data, most-often consisting of hundreds of millions of measurements/features along space and time, one has to confront the so-called  ``curse of dimensionality''. Problems arising from the ``curse of dimensionality'' in machine learning can manifest in a variety of ways, such as the excessive sparsity of data, the multiple comparison problem, the multi-collinearity of data, or even the poor generalization of the constructed surrogate models \citep{phinyomark2017resting,altman2018curse}.
To deal with these problems, methodological advances have been made involving random field theory and strict (e.g. corrected for multiple comparison problem) hypothesis testing 
to make reliable inferences in a unified framework (e.g. the Statistical Parametrical Mapping (SPM) \citep{penny2011statistical}). To this day, this approach serves as the standard approach of processing fMRI time-series, relying mainly on a Generalized Linear Model approach to infer about the activity of certain voxels or regions of them. Despite the fact that SPM has been instrumental in analyzing task-related fMRI data, at the same time imposes critical assumptions about the structure and properties of data and, thus, setting limits to the problems that can be solved \citep{madhyastha2018current}.\par 
Another effective way of dealing with problems related to the curse of dimensionality of brain signaling activity, is to represent the dataset on a low-dimensional embedding space \citep{belkin2003laplacian,coifman2006diffusion,ansuini2019intrinsic} using non-linear manifold learning \citep{gallos2021construction}. For example, \cite{qiu2015manifold} used manifold learning for the prediction of dynamics of brain networks with aging by imposing a Log-Euclidean Riemanian manifold structure on the brain and then using a framework based on Locally Linear Embedding \citep{roweis2000nonlinear} to uncover the manifold. \cite{pospelov2021laplacian} employed both linear and non-linear dimensionality reduction techniques with the aim of discriminating  fMRI resting state recordings of subjects between acquisition and extinction of an experimental fear condition. The differences in terms of interclass-intraclass distance and the classification process took part on the low dimensional space, retaining up to 10 dimensions. Non-linear manifold learning algorithms outperformed the linear ones with Laplacian eigenmaps \citep{belkin2003laplacian} producing the best results. \cite{gao2021nonlinear} used Diffusion Maps to demonstrate that task-related fMRI data of different nature span to a similar low-dimensional embedding. In a series of works, \citep{gallos2021construction,gallos2021isomap,gallos2021relation}, we have applied linear and non-linear manifold learning techniques (such as Multi-Dimensional Scaling \citep{cox2008multidimensional}, ISOMAP \citep{tenenbaum2000global}, Diffusion Maps \citep{coifman2006diffusion,nadler2006diffusion,nadler2008diffusion}) to construct embedded functional connectivity networks from fMRI data. All the above studies provided evidence that dimensionality reduction may lead to a better classification performance, i.e., a finer detection of biomarkers. However, the above studies focused more on the clustering/classification problem than building dynamical models to predict the brain dynamics. \cite{thiem2020emergent} presented a data-driven methodology based on DMs to discover coarse variables, ANNs and Geometric Harmonics to learn the dynamic evolution laws of different versions of the Kuramoto model.\par 
Here, we propose a three-step machine learning approach for the construction of surrogate dynamical models from real task fMRI data in the physical/ambient high dimensional space, thus dealing with the curse of dimensionality by learning the dynamics of brain activity on a low-dimensional space. In particular, we first use Diffusion Maps to learn an effective non-linear low-dimensional manifold. Then we learn and predict the embedded brain dynamics based on the set of variables that span the low-dimensional manifold. For this task, we use both FNNs and the Koopman Operator framework \citep{mezic2013analysis,williams2015data,brunton2016koopman,li2017extended,bollt2018matching,dietrich2020koopman,lehmberg2021modeling}. Finally, we solve the pre-image problem, i.e. the reconstruction of the predictions in the original ambient fMRI space, using Geometric Harmonics (GA) \citep{coifman2006geometric,dsilva2013nonlinear,papaioannou2021time,evangelou2022double} for the FNN predictions, and the Koopman modes per se \citep{li2017extended,lehmberg2021modeling} for the Koopman operator predictions.\par 
For our illustrations, we used a benchmark fMRI dataset, recorded during an attention to a visual motion task. The proposed approach provides out-of-sample long-term predictions, thus reconstructing accurately in the ambient space the brain activity in response to the visual/attention stimuli in the five  brain regions that are known to be the key ones for the particular task from previous studies. Furthermore, we show that the Koopman operator numerical-analysis/based scheme provides numerical approximations with similar accuracy to the FNN-GH scheme, while bypassing both the construction of a non-linear surrogate model on the DMs, and the use of GH to extrapolate predictions in the ambient fMRI space.\par 

\section{Benchmark dataset and Signal extraction}
For our illustrations, we focus on a benchmark fMRI dataset of a single subject during a task that targets attention to visual motion. This particular dataset has served as a benchmark in many studies \citep{buchel1997modulation,friston2003dynamic,penny2004comparing,almpanis2020construction} and is publicly available for download from the official SPM website \href{https://www.fil.ion.ucl.ac.uk/spm/data/attention/}{https://www.fil.ion.ucl.ac.uk/spm/data/attention/}. The images are already smoothed, spatially normalised, realigned and slice-time corrected.\par
In the experiment \citep{buchel1997modulation}, subjects observed a black screen displaying a few white dots. The experimental design consisted of specific epochs where the dots behave as either static or moving. In between these epochs, there was also an intermediate epoch where only a static picture without dots was present (i.e. a fixation phase, which is treated as baseline). The subjects were guided to keep focused for any change concerning the moving dots, even, when no changes in the intermediate epochs really existed. Thus, there were four distinct experimental conditions: ``fixation'', ``non attention'' where the dots were moving, but subjects did not need to pay attention on the screen, ``attention'' where subjects needed to pay attention on existing (if any) changes regarding the movement of dots (i.e., possible changes in velocity or acceleration) and ``static'' where the dots were stationary. The dataset consists of 4 runs stuck together, where the 10 first scans of each run were discarded for the elimination of non-desirable magnetic effects. Thus, the length of the fMRI dataset is 360 scans.\par 
Here, we have first parcellated the brain into 116 brain regions of interest (ROIs) as derived by the use of Automated Anatomical Labeling (AAL) \citep{tzourio2002automated} from the fMRI data. After the formulation of the ROIs by AAL, the average of the Blood Oxygen Level Dependent (BOLD) signal from all voxels located in each region was calculated at each time point. Out of all 116 ROIs and due to the limited field of view of the fMRI acquisition, some parts of the cerebellum were excluded as there was no existing signal (i.e. there was barely any change in the BOLD signal throughout the experiment). Specifically, these regions were the right and left part of Cerebellum\_10, Cerebellum\_7b and left part of Cerebellum\_9. The remaining 111 time series were linearly detrended and standardized for further analysis.
\section{Methodology}
\begin{figure*}[hbt!]
\centering
\includegraphics[width=1\linewidth]{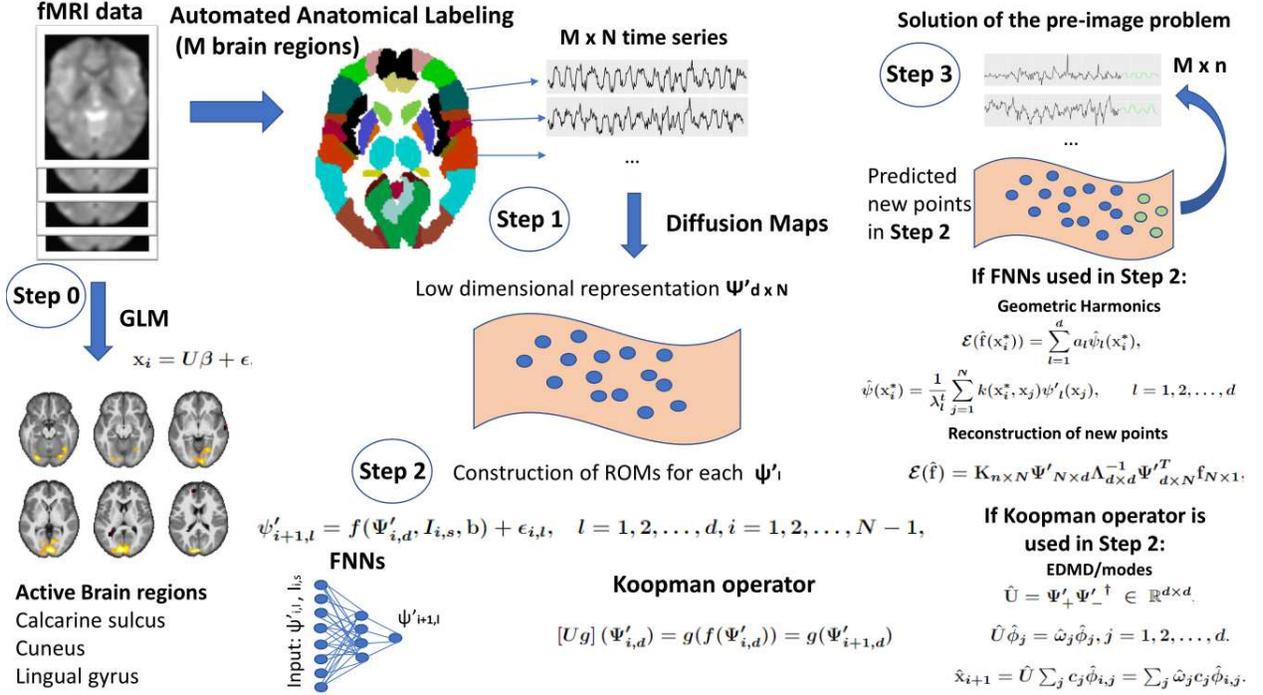}
\caption{Schematic of the proposed machine learning based methodology for modelling and predicting the dynamics of brain activity from task-dpended fMRI data.}
\label{fig:Schema}
\end{figure*}

Our data-driven machine-learning methodology for modelling grain activity deploys in three main steps. First, as a zero-step, we follow the standard procedure regarding fMRI data preprocessing with the general linear model to identify the voxels of statistically significant BOLD activity. Then, at the first step, we use parsimonious Diffusion Maps to identify a set of variables that parametrize the low-dimensional manifold where the embedded fMRI BOLD time-series evolve. Based on them, we  construct reduced-order models (ROMs), based either on FNNs or Koopman operator. Finally, at the third step, we use either Geometric Harmonics (GH) \citep{coifman2006geometric}, when we use FNNs, or Extended Dynamic Model Decomposition (EDMD) \citep{williams-2015b} when we use Koopman operator, to solve the pre-image problem, i.e. to provide out-of-sample predictions in the ambient phase space of BOLD signals. In the following sections, we describe the above steps in detail. A schematic of the three-step methodology is also shown in Fig. \ref{fig:Schema}.
\subsection{Step 0. Data Preprocessing using the General Linear Model}
For the past 20 years, the General Linear Modelling (GLM) approach has been the cornerstone of standard fMRI analysis \citep{friston1994statistical,friston1995analysis,worsley1995analysis}. Within the GLM framework, time series associated with BOLD signal at the voxel level are modeled as a weighted sum of one or more predictor variables (i.e. in relationship to the experimental conditions) including an error term for the unobserved/unexplained information in the model. The ultimate goal of the aforementioned approach is to determine whether a predictor variable does contribute (and to what extent) to the variability observed in the data (e.g., a known pattern of stimulation induced by specific experimental settings).\par 
Specifically, we can model the amplitude of the BOLD signal with fMRI in the $i$-th brain region, as a weighted sum of some known predictor variables $\bold{u}_1,...,\bold{u}_p$ each scaled by a set of parameters \boldmath$\beta$:\\
\begin{equation}
\begin{array}{l}
x_{i,1} = u_{1,1} \beta_{1} +  u_{1,2} \beta_{2} + ...  u_{1,p} \beta_{p}+ \epsilon_{i,1}\\
x_{i,2} = u_{2,1} \beta_{1} +  u_{2,2} \beta_{2} + ...  u_{2,p} \beta_{p}+ \epsilon_{i,2}\\
...\\
x_{i,N} = u_{N,1} \beta_{1} +  u_{N,2} \beta_{2} + ...  u_{N,p} \beta_{p}+ \epsilon_{i,N}
\end{array}
\end{equation}
A more compact notation leads to the following formula:
\begin{equation}
\bold{x}_i = U \bold{\beta}+ \epsilon,
\end{equation}
where $\bold{x}_{i}=[x_{i,1},x_{i,2},\dots,x_{i,N} ]^T \in \mathbb{R}^N$ is the column vector that stores the BOLD responses in all time instances, $\bold{U}=\begin{bmatrix}
u_1,1 & \dots & u_1,p\\
\cdots & \dots & \cdots\\
u_N,1 & \dots & u_N,p
\end{bmatrix} \in \mathbb{R}^{N \times p}$ is the design matrix with the column $j$-th representing the time series of the predictor variable $j$ and $\bold{\beta}=[\beta_1,\beta_2,\dots,\beta_p]^T \in \mathbb{R}^p$ is a vector of the unknown parameters that set the contribution of each of the predictors to the BOLD signal $\bold{x}_i$. Finally, $\bold{\epsilon}_i=[\epsilon_{i,1},\epsilon_{i,2},\dots,\epsilon_{i,N}]^T \in \mathbb{R}^N$ is a vector containing the corresponding (modelling) errors.\par 
In this study, the standard GLM approach was used inside the SPM framework \citep{penny2011statistical} and its toolbox extension AAL \citep{tzourio2002automated}. Out of the 3 available methods for automatic anatomical labelling of the global functional activation map, we used the cluster labeling. This is a rational choice since we wanted to localize the activated clusters in the level of anatomical regions. The final results contain which percentage of the activated cluster belongs to a specific anatomical region, and also which percentage of the region is part of the activated cluster \citep{tzourio2002automated}.\\
Here, the experimental variables $\bold{x}_i$  are the following: ``fixation'', ``static'', ``attention'', ``non-attention''. We considered as activated clusters, 10 or more nearby voxels that passed the contrast of ``attention'', ``non attention'' and ``static'' vs  ``fixation'' which is here treated as a low level baseline \citep{buchel1998functional}. We considered two different thresholds, one of $p < 0.05$ corrected for a family wise error (FWE) and another more liberal threshold of $p < 0.01$ (uncorrected). The analysis was restricted only on voxels included in the predefined regions of AAL (eg. by using the atlas as an inclusive mask). In this way, we wanted to determine the anatomical regions that are actually ``active'' throughout the experiment in respect with any of the experimental variables other than fixation (e.g., when either when subject is paying or not paying attention to moving or static dots). 

\subsection{Step 1: Identification of an Embedded Manifold in the fMRI BOLD signals via Parsimonious Diffusion Maps}
\label{sec:dmaps}
Diffusion maps is a non-linear dimensionality reduction/ manifold learning algorithm introduced by \cite{coifman2006diffusion} that exploits the inherent relationship among heat diffusion and random walk Markov chain. The main idea is that jumping from one data point to another in a random walk regime, it is more likely to jump to nearby points than points that are further away. The algorithm produces a mapping of the dataset in the Euclidean space, whose coordinates are computed utilizing the eigenvalues and the eigenvectors of a diffusion operator on the data. An embedding in the low-dimensional diffusion maps space is obtained by the projections on the eigenvectors of a normalized graph Laplacian \citep{belkin2003laplacian}. Specifically, given a dataset of $N$ points, $\bold{X}=[\bold{x}_1 \quad \bold{x}_2 \quad \dots \quad  \bold{x}_N]$, $\bold{x}_{i} \in \bold{R}^M$ is the column vector containing the BOLD signal in $N$ instances of the $i$-th brain region, $M$ corresponds to the number of brain regions, we construct an affinity matrix $W \in \mathbb{R}^{N\times N}$, defined by a kernel, any square integrable function, $w \in L^2: \mathbf{X}\times \mathbf{X}\rightarrow \mathbb{R}+^{0}$, between all pairs of points $\bold{x}_{i}$ and $\bold{x}_{j}$. The affinity matrix is most often computed though the use of the so-called Gaussian heat kernel with elements:
\begin{equation}
\label{eq:gaussian-heat-kernel}
w_{i,j}= \operatorname{exp}\left(- \frac{||\bold{x}_{i}-\bold{x}_{j}||_2^2}{2\sigma} \right),
\end{equation}
where $\vert\vert . \vert\vert_2$ computes the distance between two points using the Euclidean $L_2$ norm (any metric could be used) between points, and $\sigma$ serves as a scale parameter of the kernel. A Gaussian heat kernel in Eq. \ref{eq:gaussian-heat-kernel} produces both symmetric and positive semi-definite affinity square matrix $W$. The property is fundamental and the one that allows for the interpretation of weights as scaled probabilities. We thus formulate the diagonal normalization matrix $\bold{K}$ with elements
\begin{equation}
k_{ii}=\sum_{j=1}^{N} w_{ij}.
\end{equation}
Here, a family of anisotropic diffusions can be parameterized by adding a parameter $\alpha$ which controls the amount of influence of the density in the data. This constant can be used to normalize the affinity matrix $\bold{W}$ as:
\begin{equation}
\label{eq: normalize1}
\bold{ \Tilde{W}}= \bold{K}^{-\alpha}W\bold{K}^{-\alpha}
\end{equation}
For $\alpha = 0$, we get the normalized graph Laplacian,  $\alpha = 0.5$, we get the Fokker–Planck diffusion, and for $\alpha = 1$, we get the Laplace–Beltrami operator \citep{thiem2020emergent}.\par 
Next, we renormalize the matrix $\bold{\Tilde{W}}$ with the diagonal matrix $\Tilde{K}$, with elements $\Tilde{k}_{ii}=\sum_{j=1}^{N} \Tilde{w}_{ij}$. $\bold{P} \in \mathbb{R}^{N\times N}$ is a Markovian row stochastic matrix which is also called diffusion matrix given by: 
\begin{equation}
\label{eq:normalize2}
\bold{P}=\bold{\Tilde{K}}^{-1}\bold{\Tilde{W}}.
\end{equation}
In this stage, the elements of the diffusion matrix can be seen as the scaled probabilities of one point jumping to another in a random walk sense. Consequently, taking the power $t$ of the diffusion matrix $\bold{P}$ is essentially identical of observing $t$ steps forward of a Markov chain process on the data points. The element $\bold{P}^t(\bold{x}_i,\bold{x}_j)$ denotes the transition probability of jumping from point $\bold{x}_i$ to point $\bold{x}_j$ in $t$ steps.\par 
Applying singular value decomposition (SVD) on  $\bold{P}$, we get: 
\begin{equation}
\bold{P}=\bold{\Psi}\bold{\ \Lambda} \bold{\Psi}^{T},
\label{eq:decomp}
\end{equation}
with $\bold{\Lambda}$ being the diagonal matrix that stores the $M$ eigenvalues and $\bold{\Psi}\in \mathbb{R}^{N \times N}$ the eigenvectors of $\bold{P}$.\par 
These eigenvectors are the discrete approximations of the eigenfunctions of the corresponding to the parameter $\alpha$ (see above) kernel operator on the manifold. Thus, such eigenfunctions provide an orthonormal basis of the $L^2$ space, i.e. of the space of square-integrable functions \citep{coifman2006diffusion}.\par 
The eigenvalues that correspond to the eigenvectors in descending order are $\lambda_{0}=1 \geq \lambda_{1} \geq \lambda_{2}$ ... $\geq \lambda_{k}$ with the first being the trivial eigenvalue that is equal to 1 ($\bold{P}$ is a Markovian matrix) and $k$ being the target embedding dimension. The embedding is obtained by mapping the points from the original/ ambient space to the diffusion maps space by preserving the diffusion distance among all data points in the dataset \citep{nadler2006diffusion}. Since the diffusion distance takes into account all possible paths between data points, it is thus more robust against noise perturbations \citep{coifman2006diffusion} (comparing for example with the geodesic distance).\par 
The coordinates of the nonlinear projection of a point $\bold{x}_{i}$ in the $k$ dimensional space spanned by the first $k$ DMs $\{\bold{\psi}_1,\bold{\psi}_2,\dots, \bold{\psi}_k\}$ read:
\begin{equation}
\mathcal{R}(\bold{x}_{i})\equiv \bold{y}_{i} =[y_{i,1},y_{i,2},\dots,y_{i,k}]^T=[\lambda^t_{1}\psi_{i,1},\lambda^t_{2}\psi_{i,2},\dots,\lambda^t_{k}\psi_{i,k}]^T, i=1,2,\dots,N.
\label{eq:DMcoords}
\end{equation}
where, $\psi_{i,l}$ is the $i$-th element of $\bold{\psi}_l, l=1,2,\dots,k$.
For all practical purposes, the embedding dimension is determined by the spectral gap in the eigenvalues of the final decomposition. A numerical gap (or a sharp decrease) between the first few eigenvalues and the rest of the eigenspectrum would indicate that a few eigenvectors could be adequate for the approximation of the diffusion distance between all pairs of points \citep{coifman2006diffusion}.
A compact version of the algorithm is the following:
\begin{enumerate}
    \item \textbf{Input:} linearly detrended data set $\bold{X}=[\bold{x}_1 \quad \bold{x}_2 \quad \dots \quad  \bold{x}_N]$, $\bold{x}_{i} \in \mathbb{R}^M$.
    \item Compute the Gaussian heat kernel with elements $w_{i,j}= \operatorname{exp}\left(- \frac{||\bold{x}_{i}-\bold{x}_{j}||_2^2}{2\sigma} \right), i,j=1,2,\dots,N$.
    \item Apply the normalization on $\bold{W}$ (Eq. \ref{eq: normalize1}) and renormalization (Eq. \ref{eq:normalize2}) to compute the diffusion matrix $\bold{P}$; 
    \item Apply SVD on matrix $\bold{P}$ (Eq. \ref{eq:decomp}) and keep the first $k$ eigenvectors corresponding to the $k$ largest eigenvalues.
    \item \textbf{Output:} $\bold{\Psi}_{k}=[\bold{\psi}_1 \quad  \bold{\psi}_2 \quad  \dots \bold{\psi}_k] \in \mathbb{R}^{N\times k}$, $\bold{y}_{i} =[y_{i,1},y_{i,2},\dots,y_{i,k}]^T$, $i=1,2,\dots,N$ (see Eq.(\ref{eq:DMcoords}).
\end{enumerate}
Next, from the initial set of eigenvectors, we retain only the $d$ parsimonious eigendimensions as proposed in \citep{dsilva2018parsimonious}. Specifically, using a local linear function,
\begin{equation}
\psi_{i,l}\approx c_{i,l}+\bold{\beta}_{i,l}^{T} \Psi_{i,l-1},\quad \Psi_{i,l-1}=[\psi_{i,1}, ... ,\psi_{i,l-1}]^T \in \mathbb{R}^{l-1},
\end{equation}
$\alpha_{i,l} \in \mathbb{R}, \bold{\beta}_{i,l} \in \mathbb{R}^{l-1}$ we approximate each point of $\psi_{i,l}$ based on the remaining data points $\psi_{i,l-1}$. 
This leads to the following optimization problem \citep{thiem2020emergent}:
\begin{equation}
\label{optprob}
\argmin_{c,\bold{\beta}} \sum\limits_{ \ i \neq j } K\Big(\Psi_{i,l-1},\Psi_{j,l-1}\Big)\Big(\psi_{j,l}-(c+\bold{\beta}^T \Psi_{j,l-1})\Big)^2,
\end{equation}
where  $K$ is the Gaussian kernel. Finally, the computation of the normalized error $er_k$ is done through the use of leave-one out cross validation for each of the local linear fit. The normalized error is defined as:
\begin{equation}
\label{normerr}
er_l=\sqrt{\frac{\sum_{i=1}^{N}\Big(\psi_{i,l}-(\alpha_{i,l}+\bold{\beta}_{i,l}^T\bold{\Psi}_{i,l-1})\Big)}{\sum_{i=1}^{N}(\psi_{i,l})^2}}.
\end{equation}
A small or negligible error $er_l$ denotes that $\bold{\psi}_l$ can be actually predicted from the remaining eigenvectors $\bold{\psi}_1, \bold{\psi}_2,...,\bold{
\psi}_{l-1}$ and thus is a repeated eigendirection (i.e., $\bold{\psi}_l$ is considered a harmonic of the previous eigenmodes). Therefore, only the eigenvectors that exhibit a large $er_l$ are selected in a way of seeking the most parsimonious representation. More information regarding the choice of the most parsimonious eigenvectors can be found in \citep{dsilva2018parsimonious}, \citep{lee2020coarse} and \citep{galaris2022numerical}.\par 
A compact pseudocode for identifying parsimonious eigenvectors is the following:
\begin{enumerate}
    \item \textbf{Input:} Set of $\bold{\Psi}_k \in \mathbb{R}^{N\times k}$, $k<<M$ corresponding to the largest eigenvalues from the application of Diffusion Maps and $d$ the number of parsimonious eigenvectors to retain.
    \item Solve the optimization problem given by Eq.(\ref{optprob}) and calculate the normalized error $er_l$ given by Eq.(\ref{normerr}).
    \item \textbf{Output:} $\bold{\Psi'}_d\in \mathbb{R}^{N\times d}$, with columns the parsimonious DMs $\mathbf{\psi}'_j \in \mathbb{R}^N, j=1,2,\dots,d$ corresponding to the $d$ largest $er_l$s.
\end{enumerate}

\subsection{Step 2. Construction of Diffusion Maps-based Reduced order models}
\label{sec:reduced}
After the reduction of space dimension to a set of $d$ of parsimonious DMs components, we proceed to the construction of reduced order (regression) models (ROMs).\par 
Here, we have used two approaches for the construction of ROMs and lifting back to the original ambient space: (a) via  FNNs and Geometric Harmonics, and, (b) via the Koopman modes/EDMD.
In both cases, the general form of a ROM reads:
\begin{equation}
\label{roms}
\psi'_{i+1,l}= f(\mathbf{\Psi}'_{i,d},I_{i,s},\bold{b})+\epsilon_{i,l}, \quad l=1,2,\dots,d, i=1,2,\dots,N-1,
\end{equation}
where $\mathbf{\Psi}'_{i,d}=[\mathbf{\psi}'_{i,1},\dots,\mathbf{\psi}'_{i,l},\dots,\mathbf{\psi}'_{i,d}]^T \in \mathbb{R}^d$, is the vector containing the $i$-th elements of the  $d$ parsimonious DMS ($\bold{\psi'}_{i,l}$ denotes the $i$-th element of the $l$-th parsimonious DM), $b$ denotes a set of parameters of the model, $\bold{I}_{i,s}$ the external stimuli $s$ indicating specific experimental conditions at the time instant $i$, and $\epsilon_{i,l}$ is the modelling error. We note that after training, predictions were done iteratively: we set the initial conditions at the last point of the training set, and then iterate the ROMs to produce the long-term predictions.\par 
In order to compare the efficiency of the two approaches, we applied also a simple naive random walk (NRW) model in a direct way for one-step ahead predictions: predictions at time $i+1$ are equal to the last observed values, i.e., $\psi_{i+1,l}=\psi_{i,l}$. \par 
Thus, in contrast to the predictions via the embedded ROMs that are produced iteratively, i.e. without any knowledge of the previous values, but the initial conditions, the predictions made via the NRW model are based on the knowledge of the values of the previous points of the test set.

\subsubsection{Reduced order models with FNNs}
As a first reduced order modelling approach on the discovered manifold $\mathcal{M}$, we used $d$ FNNs, for each one of the embedded coordinates $\bold{\psi}'_l, l=1,2,\dots,d$, with one hidden layer, $H$ units and a linear output layer with $n$ units that can be written compactly as
\begin{equation}
\begin{aligned}
    {\psi'}_{i+1,l}=
     \bold{w}^T_{l,o} \bold{S}( \bold{W}^T_{1,l} \mathbf{\Psi}'_{i,d} + \bold{b}_{1,l})+b_{o,l}, \quad  l=1,2,\dots,d.
\end{aligned}
\label{eq:FNN}
\end{equation}
$\mathbf{S(\cdot)}$ is the activation function (based on the above formulation it is assumed to be the same for all networks and nodes in the hidden layer), $\mathbf{b}_{1,l} \in \mathbb{R}^{H}$ is the vector containing the biases of the nodes of the first hidden layer, $\bold{W}_{1,l}\in\mathbb{R}^{d\times H}$ is the matrix containing the weights between the input and the first hidden layer, $\bold{w}_{o,l} \in \mathbb{R}^{H}$ denotes the vector containing the weights between the hidden layer and the output layer, and finally $b_{o,l}$ is the bias of the output.\par 
The total number of input units were matched to the embedding coordinates of the manifold plus the external stimuli $\bold{I}$ (i.e., $d$ embedding coordinates plus 4 external stimuli: attention, non attention, static and fixation) while there was one output unit. We used the logistic transfer function $S(x)=\frac{1}{1+e^{-x}}$ as the activation function for all neurons in the hidden layer. A learning rate decay parameter was also used as a regularization parameter to prevent over-fitting and improve generalization \citep{krogh1992simple} of the final model. The neural networks were trained on the first 280 time points (which accounts for roughly the 77\% of the data points) using repeated (10 times) 10-fold cross validation approach. All parameters like different number of neurons $a$ and learning rate decay parameter $c$ were optimized via grid search inside the cross validation procedure.
Finally, we optimized five final models (i.e., one model for each one of the parsimonious eigenvectors selected) to make one-step predictions in the following way: we give $\psi_{i,l}$ as input to the FNN in order to make prediction on $\psi_{i+1,l}$.
Using the best candidate model for each one of the parsimonious components $\bold{\psi}_{l}$, we predicted iteratively the left out/ unseen test points (i.e., the next 80 unseen time points).  

\subsubsection{Reduced order models with the Koopman operator}
\label{sec:koopman}
In the Koopman operator framework, predictions are performed in the function space over the data manifold, and so state prediction turns into coordinate function prediction~\citep{mezic-2005,budisic-2012,dietrich-2020a}. The defining property in this framework is that the operator always acts linearly on its domain, which makes it amenable to spectral analysis and approximation techniques from numerical analysis/ linear algebra.
Given the flow map $f:\mathbb{R}^d\rightarrow \mathbb{R}^d$ of the embedded dynamics on the manifold, such that $\mathbf{\Psi}'_{i+1,d}=f
(\mathbf{\Psi}'_{i,d})$, the Koopman operator $U:\mathcal{F}\to\mathcal{F}$ acts on observables $g\in\mathcal{F}$ by composition with $f$, such that
\begin{equation}
\left[U g\right](\mathbf{\Psi}'_{i,d})= g(f(\mathbf{\Psi}'_{i,d}))=g(\mathbf{\Psi}'_{i+1,d}).
\end{equation}
If a function $\phi_j$ is an eigenfunction of $U$ with eigenvalue $\omega_j$, then $\left[U\phi_j\right](\mathbf{\Psi}'_{i,d})=\phi_j(\mathbf{\Psi}'_{i+1,d})=\omega_j \phi_j(\mathbf{\Psi}'_{i,d}).$ For many systems, eigenfunctions of the Koopman operator span a large subspace in $\mathcal{F}$, so that many observables $g\in\mathcal{F}$ can be written as a linear combination of eigenfunctions, $g=\sum_j c_j\phi_j$, with coefficients $c_j\in\mathbb{C}$. These coefficients are often called ``Koopman modes'' of the observable $g$.\par 
The most prevalent numerical algorithm to approximate the Koopman operator is dynamic mode decomposition  (DMD)~\citep{schmid-2010,schmid-2022}. It approximates a linear map between the given observables of a system and their future state, assuming that the observables span a function subspace invariant to the dynamics. In our case, we apply DMD to five of the embedding diffusion map coordinates $\bold{\Psi'}_d$ (this can be interpreted as a vector-valued observable, i.e., the Koopman operator is acting on five coordinates simultaneously~\citep{budisic-2012}). 

\subsection{Step 3. Predictions in the Ambient fMRI Space: Solution of the Pre-Image Problem}
 \label{sec:lift}
 The final step is to ``lift''predictions back to the high dimensional ambient space (here a 111-dimensional space of brain regions), i.e. solve the so-called ``pre-image'' problem.\par 
For linear manifold learning methods such as the Principal Component Analysis (PCA), this is trivial. However, for non-linear manifold learning algorithms such as DMs, there is no explicit inverse mapping and no unique solution (see \cite{chiavazzo2014reduced,papaioannou2021time,evangelou2022double}) and one has to construct a lifting operator $\mathcal{L}=\mathcal{R}(\bold{X}) \rightarrow \bold{X}$ for new unseen samples on the manifold $\bold{y}_{*} \notin \mathcal{R}(\bold{X})$.\\
The problem regarding an unseen data point is often referred to as the ``out-of-sample extension'' problem. While this problem traditionally addresses the direct problem (i.e., the extension of a function $g$ on the ambient space for new unseen data points $\bold{x}_{*}$ $\notin$ $\bold{X}$), here, we are interested in the inverse problem (i.e., the lifting of predictions made on the manifold back to the ambient space). Here, as discussed, the out-of-sample solution of the pre-image problem is solved using Geometric Harmonics when using FNNs for ROMs, or with Extended Dynamics Mode Decomposition when using the Koopman operator. Below, we present both methods. 

\subsubsection{Geometric Harmonics}
Typically, the out-of-sample extension problem is solved through the Nystr\"om extension methodology \citep{nystrom1929praktische} derived from the solution of the Fredholm equation of the second kind: 
\begin{equation}
\label{fredholm}
f(t)=g(t)+\mu \int_{a}^{b} k(t,s)f(s)ds,
\end{equation}
where $k(t,s)$ and $g(t)$ are known functions while $f(t)$ is the unknown. The Nystr\"om method approximates the integral 
\begin{equation}
\label{eqapprox}
\int_{a}^{b} y(s)ds\approx \sum_{j-1}^Nw_{j}y(s_j),
\end{equation}
where $N$ are the collocation points and $w$ the weights determined by, for example, the Gauss-Jacobi quadratic rule. Using Eq.\ref{eqapprox} in \ref{fredholm} and evaluating $f$ and $g$ at the collocation points $N$, we can make the following approximation:
\begin{equation}
(\bold{I}-\mu \Tilde{K})\hat{\bold{f}}=\bold{g},
\end{equation}
where $\Tilde{K}$=$\Tilde{k}_{ij}=k(s_i,s_j)w_j$. Solution to the homogenous Fredholm problem (where $\bold{g=0}$) is provided through the solution of the eigenvalue problem:
\begin{equation}
\Tilde{K}\hat{\bold{f}}=\frac{1}{\mu}\hat{\bold{f}},
\end{equation}
i.e.,
\begin{equation}
\sum_{j=1}^Nw_{j}y(s_i,s_j)\hat{f_j}= \frac{1}{\mu}\hat{f_i}, \qquad  i=1,2,\dots,N
\end{equation}
where $\hat{f_i}=\hat{f}(s_i)$ the $i$-th component of $\hat{f}$. The Nystr\"om extension in the full domain for  $N$ number of collocation points at a point $x$ is given by:
\begin{equation}
\label{map}
\mathcal{E}(f(x))=\hat{f}(x)=\mu\sum_{j-1}^Nw_{j}k(x,s_j)\hat{f_j}.
\end{equation}
The above Eq.(\ref{map}) provides a map from high dimensional to the reduced order space and back.\par 
Since the eigenvectors of the DMs form a basis for the embedded manifold, the extension is formulated \citep{coifman2006geometric} by the expansion of $\bold{f}(\bold{x}_i)$ using the first $d$ parsimonious eigenvectors of the diffusion matrix $\bold{P}^{t}$:
\begin{equation}
\hat{f}(\bold{x}_i)=\sum_{l=1}^da_{l}{\bold{\psi'}}_l(\bold{x}_i),\qquad i=1,2,\dots,N,
\end{equation}
where $a_l=\langle {\bold{\psi'}}_l,\bold{f}(\bold{x}_i)\rangle$, ($\langle \cdot,\cdot \rangle$ denotes inner product), $\bold{f} \in \mathbb{R}^N$ is the vector containing the values of the function $f$ at the $N$ points. Thus, the Nystr\"om extension of $\bold{f}$ at a new unseen point $\bold{x}^{*}_i$ using the same coefficients $a_k$, reads:
\begin{equation}
\mathcal{E}(\hat{\bold{f}}(\bold{x}^{*}_i))=\sum_{l=1}^da_{l}\hat{\bold{\psi}}_l(\bold{x}^{*}_i),
\label{eq:Nystrom}
\end{equation}
where,
\begin{equation}
\hat{\bold{\psi}}(\bold{x}^{*}_i)=\frac{1}{\lambda_l^t} \sum_{j=1}^Nk(\bold{x}^{*}_i,\bold{x}_j)\bold{\psi'}_l(\bold{x}_j), \qquad  l=1,2,\dots,d,
\end{equation}
are the geometric harmonics extension of each one of the parsimonious eigenvectors to the unseen data. Thus, for a new set $\bold{X}^{*}$ of $n$ points we get:
\begin{equation}
\mathcal{E}(\hat{\bold{f}})=\bold{K}_{n\times N}\bold{\Psi'}_{N\times d}\bold{\Lambda}^{-1}_{d\times d}\bold{\Psi'}^{T}_{d\times N}\bold{f}_{N \times 1},
\label{eq:reconstruct}
\end{equation}
where $\bold{K}$ is the $n\times N$ kernel matrix, $\Lambda$ is a $d\times d $ diagonal matrix that stores the eigenvalues of the eigenvectors.\par 
Here, for the solution of the pre-image problem, we used the approach of ``double Diffusion maps''  (\citep{evangelou2022double,papaioannou2021time,patsatzis2023data}). A pseudocode for the solution of the pre-image problem is given below:
\begin{enumerate}
    \item \textbf{Input:} $\bold{\Psi'}_d\in \mathbb{R}^{N\times d}$ matrix whose columns are the $d$ parsimonious diffusion maps (from Eq.(\ref{sec:dmaps})).
    \item For each new point $\{\bold{y}^{*}_i\}^{n}_{i=1}$ in the DMs embedded space computed using the Nystr\"om extension (Eq.(\ref{eq:Nystrom}), compute the kernel matrix $\bold{K}\in \mathbb R^{n\times N}$ with elements $k_{i,j}= \operatorname{exp}\left(- \frac{||\bold{y}^{*}_{i}-\bold{y}_{j}||_2^2}{2\sigma} \right)$, $i=1,2,\dots,n$, \quad $j=1,2,\dots,N$. 
    \item \textbf{Output:} Reconstructed points in the ambient space using Eq.(\ref{eq:reconstruct}).
\end{enumerate}

\subsubsection{Koopman Operator and Extended Dynamic Mode Decomposition (EDMD)}\label{sec:EDMD}
When using the Koopman operator approach, the DMs coordinates considered as a truncated basis of a function space over the original measurements, can be used effectively similarly as when using Dynamic Mode Decomposition (DMD) (see~\citep{williams-2015b}), but on the original measurements in the ambient fMRI space.\par  
Thus, unlike standard DMD, in our case, the Koopman modes are computed for the original coordinates, not the DM eigenfunctions, so that we obtain a linear map from Koopman eigenfunctions to original measurements. A compact version of the EDMD algorithm we use is as follows:
\begin{enumerate}
    \item \textbf{Input:} Diffusion maps coordinates $ \mathbf{\Psi}'_{i,d}\in \mathbb{R}^{N\times d}$ at time steps $i=1,\dots,N$ from Section \ref{sec:dmaps}.
    \item Approximate the Koopman matrix ${\hat{\mathbf{U}}=\mathbf{\Psi}'_{+} {\mathbf{\Psi}'_{-}}^{\dagger}}  \in \mathbb{R}^{d\times d}$, where $\dagger$ denotes the pseudo-inverse and ${\mathbf{\Psi}'_{-} = \left[\mathbf{\Psi}'_{1,d}, \mathbf{\Psi}'_{2,d}, \ldots, \mathbf{\Psi}'_{N-1,d} \right]}\in \mathbb{R}^{d\times (N-1)}$ and ${\mathbf{\Psi}'_{+} = \left[\mathbf{\Psi}'_{2,d}, \mathbf{\Psi}'_{3,d}, \ldots, \mathbf{\Psi}'_{N,d} \right]}\in \mathbb{R}^{d \times (N-1)}$ are time shifted data matrices with $\mathbf{\Psi}'_{i,d}\in\mathbb{R}^{d} $ as columns.
    \item Compute all eigenvectors $\hat{\mathbf{\phi}}_j \in \mathbb{R}^d$ and eigenvalues $\hat{\omega}_j$ of $\hat{\mathbf{U}}$ by solving $\hat{U}\hat{\phi}_j=\hat{\omega}_j \hat{\phi}_j, j=1,2,\dots,d$.
    \item Obtain the Koopman modes $c_j\in\mathbb{R}^M$ for the original measurements $\mathbf{x}_i\in\mathbb{R}^M$, by solving the linear systems $\mathbf{x}_i=\sum_j c_j \hat{\mathbf{\phi}}_{i,j}$ with data from all available time steps $i=1,\dots,N$ simultaneously.
    \item \textbf{Output:} Eigenvectors $\hat{\mathbf{\phi}}_j$, eigenvalues $\hat{\omega}_j$, and Koopman modes $c_j$ for the original measurements.
    \item \textbf{Prediction:} To obtain an approximation of $\mathbf{x}_{i+1}$, we evaluate $\hat{\mathbf{x}}_{i+1}=\hat{U}\sum_j c_j\hat{\mathbf{\phi}}_{i,j}=\sum_j \hat{\omega}_j c_j \hat{\mathbf{\phi}}_{i,j}$.
\end{enumerate}
\section{Results}
For our computations, we have used the Python packages Datafold \citep{lehmberg2020datafold} with scikit-learn  \citep{pedregosa2011scikit}. For the implementation of the FNNs we utilized the R package ``nnet'' \citep{ripley2016package}.
\subsection{Data Preprocessing}
\begin{figure*}[hbt!]
\centering
\includegraphics[width=0.9\linewidth]{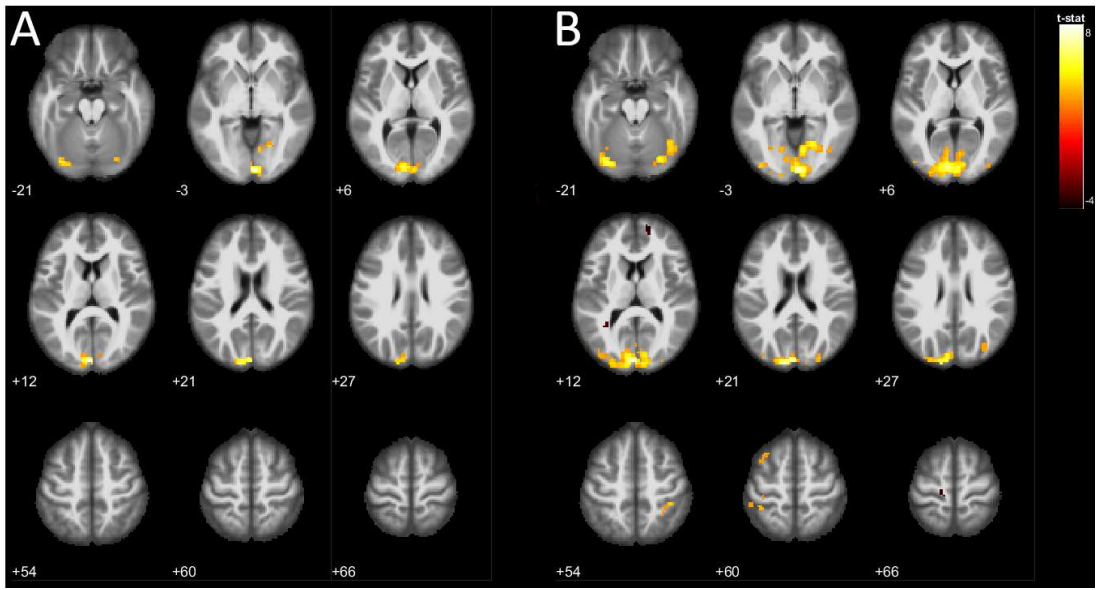}
\caption{Thresholded image of the contrast ``attention'' + ``non-attention'' + ``stationary'' vs ``fixation'' is presented as an overlay on the reference structural MRI (T1) image.A) Level of significance p<0.05 (FWE corrected), B) $p<0.001$ (uncorrected).}
\label{fig:figure1}
\end{figure*}

Our analysis starts with the data-preprocessing using GLM. In Figure \ref{fig:figure1}, the thresholded image for the activation during the contrast ``attention'' + ``non-attention'' + ``stationary'' vs ``fixation'' is presented as an overlay on the brain extracted structural T1 MNI image (2mm space). We considered activations of clusters of 10 or more voxels with two different levels of significance: A) ${p < 0.05}$ with Family Wise Error correction (FWE) and B) ${p < 0.001}$ (uncorrected). In Table \ref{tab:table1}, the results considering the more stringent threshold of $p < 0.05$ (FWE corrected) are presented in detail. We report the local maximum of each cluster, its size in number of voxels, the brain regions that contribute to each cluster and finally the percentage of the contribution. Similarly, in Table \ref{tab:table2} the results concerning the activations of clusters using a more liberal threshold of $p < 0.001$ (uncorrected) are reported. The brain regions that showed activations were mainly parts of the Occipital Lobe such as the Calcarine, Cuneus, Occipital gyrus, the Fusiform gyrus and parts of the Cerebellum. When a more liberal threshold was applied ($p < 0.001$), the clusters would naturally get larger. Consequently, two more regions appeared to have clusters of activation, namely, the Postcentral gyrus and the inferior parietal cortex. Focusing on these ``active'' regions, we then proceeded to the prediction of their behaviour through macroscopic variables extracted by the Diffusion Maps.  
\begin{table}[H]
\centering
\caption{Clusters of 10 or more activated voxels for the contrast ``attention'' + ``non-attention'' + ``stationary'' vs ``fixation'' using SPM. The level of significance was set to $p < 0.05$ (FWE corrected). For each cluster, the coordinates of the local maximum (in standard MNI coordinates) is shown along with the total number of voxels in the cluster and the brain regions that these voxels are found. Finally, the percentage of a region's voxels in each cluster is also presented.}
\begin{tabular}{llll}
\multicolumn{4}{c}{Level of significance p \textless{} 0.05 (FWE corrected)}                 \\ \hline
Local maximum x,y,z mm & Number of voxels in the cluster & Region Label      & \% in the cluster \\\hline
3, -93, 0              & 191                             & Calcarine L       & 11.96             \\
                       & 191                             & Cuneus L          & 10.17             \\
                       & 191                             & Occipital Sup L   & 8.89              \\
                       & 191                             & Calcarine R       & 4.72              \\
                       & 191                             & Occipital Mid L   & 0.21              \\
                       & 191                             & Lingual R         & 0.15              \\\hline
-33, -84, -21          & 28                              & Cerebellum Crus1 L & 39.29             \\
                       & 28                              & Lingual L         & 39.29             \\
                       & 28                              & Fusiform L        & 17.86             \\
                       & 28                              & Occipital Inf. L  & 3.57              \\\hline
21, -66, -3            & 20                              & Lingual R         & 100               \\\hline
27, -81, -18           & 11                              & Cerebellum 6 R     & 45.45             \\
                       & 11                              & Cerebellum Crus1 R & 27.27             \\
                       & 11                              & Fusiform R        & 18.18             \\
                       & 11                              & Lingual R         & 9.09
\end{tabular}
\label{tab:table1}
\end{table}
\begin{table}[H]
\centering
\caption{Clusters of 10 or more activated voxels for the contrast ``attention'' + ``non-attention'' + ``stationary'' vs ``fixation'' using SPM. The level of significance was set to p < 0.001 (uncorrected). For each cluster, the coordinates of the local maximum (in standard MNI coordinates) is shown, along with the total number of voxels in the cluster and the brain regions that these voxels are found. Finally, the percentage of a region's voxels in each cluster is also presented.}
\begin{tabular}{llll}
\multicolumn{4}{c}{Level of significance p \textless{} 0.001 (uncorrected)}                          \\ \hline
Local maximum (x,y,z mm) & Number of voxels in the cluster & Region Label      & \% in the cluster \\\hline
3, -93, 0                & 917                             & Calcarine L       & 20.17             \\
                         & 917                             & Occipital Mid L   & 13.52             \\
                         & 917                             & Lingual R         & 12.65             \\
                         & 917                             & Occipital Sup L   & 11.01             \\
                         & 917                             & Calcarine R       & 10.58             \\
                         & 917                             & Fusiform R        & 10.25             \\
                         & 917                             & Cuneus L          & 7.96              \\
                         & 917                             & Lingual L         & 6                 \\
                         & 917                             & Cerebellum Crus1 R & 2.94              \\
                         & 917                             & Cerebellum 6 R     & 1.74              \\
                         & 917                             & Cuneus R          & 1.64              \\
                         & 917                             & Occipital Inf. L  & 1.2               \\
                         & 917                             & Occipital Sup R   & 0.22              \\
                         & 917                             & Occipital Inf. R  & 0.11              \\\hline
-33, -84, -21            & 68                              & Cerebellum Crus1 L & 44.12             \\
                         & 68                              & Lingual L         & 25                \\
                         & 68                              & Fusiform L        & 19.12             \\
                         & 68                              & Occipital Inf. L  & 11.76             \\\hline
-39, -90, -3             & 10                              & Occipital Mid. L  & 90                \\
                         & 10                              & Occipital Inf. L  & 10                \\\hline
39, -39, 54              & 15                              & Parietal Inf. R   & 80                \\
                         & 15                              & Postcentral R     & 20                \\\hline
-33, -69, -18           & 35                              & Fusiform L        & 94.29             \\
                         & 35                              & Cerebellum 6 L     & 2.86              \\
                         & 35                              & Lingual L         & 2.86              \\\hline
27, -90, 20              & 11                              & Occipital Sup. R  & 72.73             \\
                         & 11                              & Occipital Mid.R   & 27.27            
\end{tabular}
\label{tab:table2}
\end{table}

\subsection{Manifold Learning with Diffusion Maps}
\begin{figure*}[hbt!]
\centering
\includegraphics[width=0.8\linewidth]{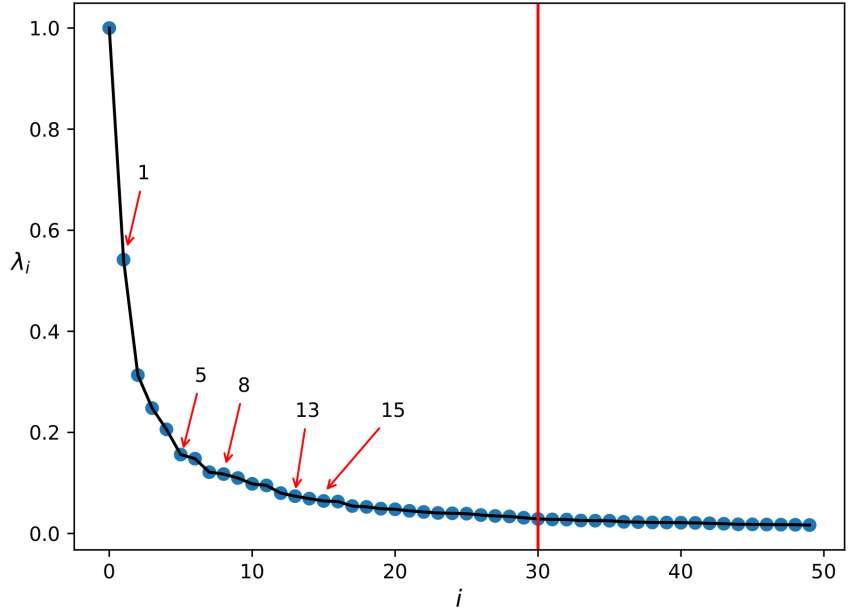}
\caption{The eigenspectrum of DMs on the 111 time series that correspond to different brain regions. The red vertical line indicates how many eigenvectors we computed. The five parsimonious eigenvectors that were finally used as macroscopic variables are also marked with red arrows.}
\label{fig:figure2}
\end{figure*}

Manifold learning with Diffusion maps was applied on the first 280 (out of the whole 360 time points) time points of the 111 linearly detrended time series. The last 80 points were left out to evaluate later the prediction error based on ROMs and the ``lifting'' process. In Figure \ref{fig:figure2}, we present the eigenspectrum of the DMs decomposition. The red vertical line shows the number of the extracted eigenvectors (more eigenvectors would not affect the final outcome). The five most parsimonious DMs components that we used for further analysis are marked with red arrows, namely $\psi_1, \psi_5, \psi_8, \psi_{13}$ and $\psi_{15}$. As it can be seen in Figure \ref{fig:figure2}, after the first 30 eigenvalues, most of the variance has already been captured. In other words, consequent eigenvectors have a very small contribution to the final embedding. The parameters we actually set in our computations are $t=0$, $\sigma=65$, $\alpha=1$, $\kappa=30$ and $d=5$. 

\subsection{Construction of Reduced Order Models and Predictions on the Manifold}
\begin{figure*}[hbt!]
\centering
\includegraphics[width=0.9\linewidth]{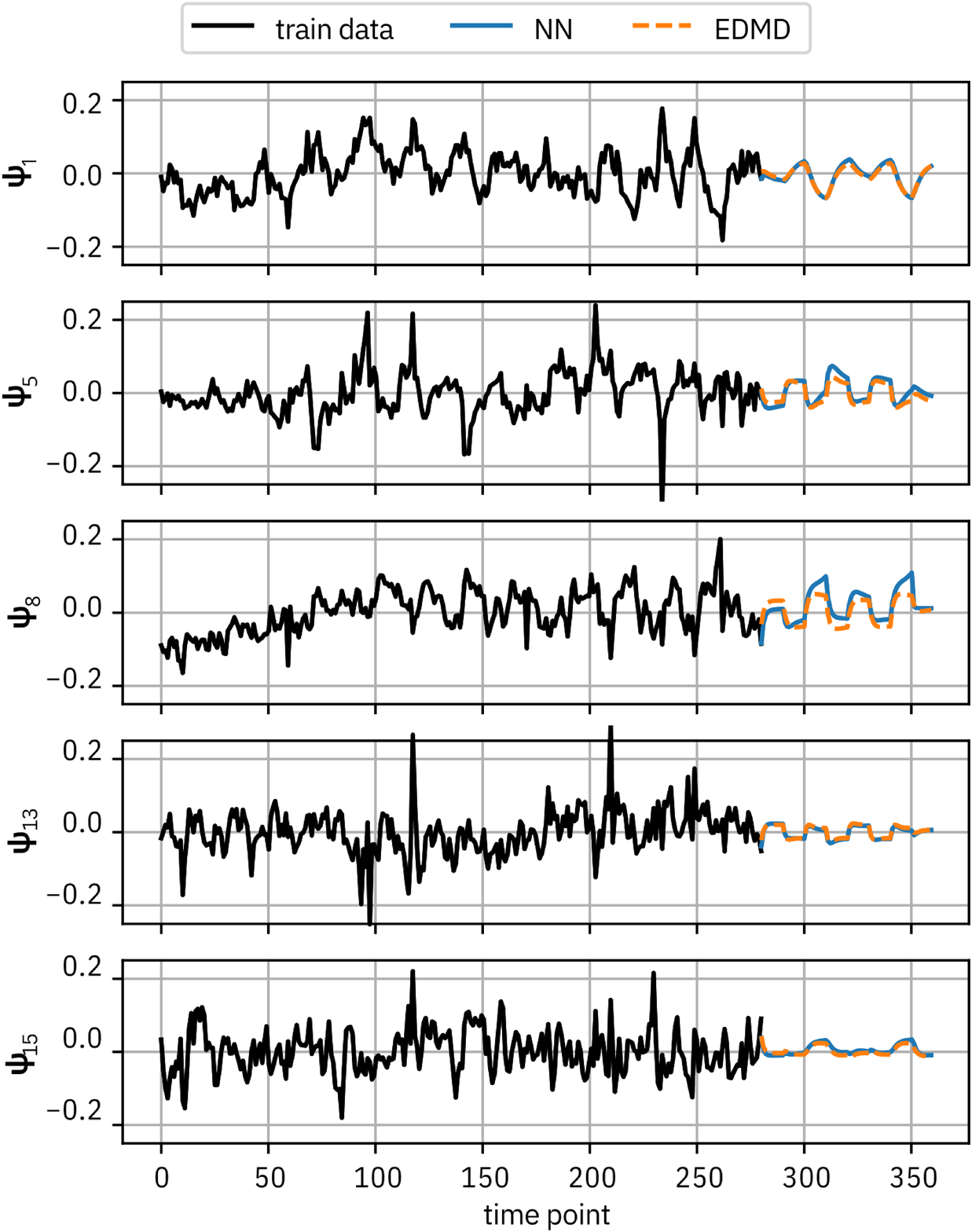}
\caption{Predictions based on the 5 parsimonious eigenvectors (namely $\mathbf{\psi}_1, \mathbf{\psi}_5, \mathbf{\psi}_8, \mathbf{\psi}_{13}$ and $\mathbf{\psi}_{15}$) when applying FNNs and the Koopman operator. Actual points for each one of the DMs are presented with black color up to the 280th point, which is the last point of the training set. The predicted values are presented with light green color up to the end of the time series (as derived iteratively by the ROMs from the point 281 to the end.).}
\label{fig:figure3}
\end{figure*}

Out-of-sample predictions on the embedded manifold were based on the 5 parsimonious DMs, namely $\mathbf{\psi}_1, \mathbf{\psi}_5, \mathbf{\psi}_8, \mathbf{\psi}_{13}$ and, $\mathbf{\psi}_{15}$ are shown in Figure \ref{fig:figure3}. Based on the these, we constructed ROMs using FNNs and the Koopman operator as described in subsection \ref{sec:reduced}. \par 
The ROMs were trained via one time step ahead predictions using the first 280 time points, which account to the 77\% of the total data points. Validation was done using 10-fold cross
validation, repeated 10 times. Thus, the performance of the ROMs was assessed by simulating the trained ROMs iteratively: setting the initial conditions at the time point 280, we iterated the ROMs to produce the next 80 points (i.e. up to the 360th time point).\par 
For the FNNs, we used the logistic function as the activation function for all neurons in the hidden layer and a learning rate decay parameter or regularization to prevent
over-fitting and improve generalization of the final model. Parameter tuning was done via grid search during the repeated cross validation process.

\subsection{Predictions in the fMRI space and comparison with the Naive Random Walk Model}
\begin{figure*}[hbt!]
\centering
\includegraphics[width=1\linewidth]{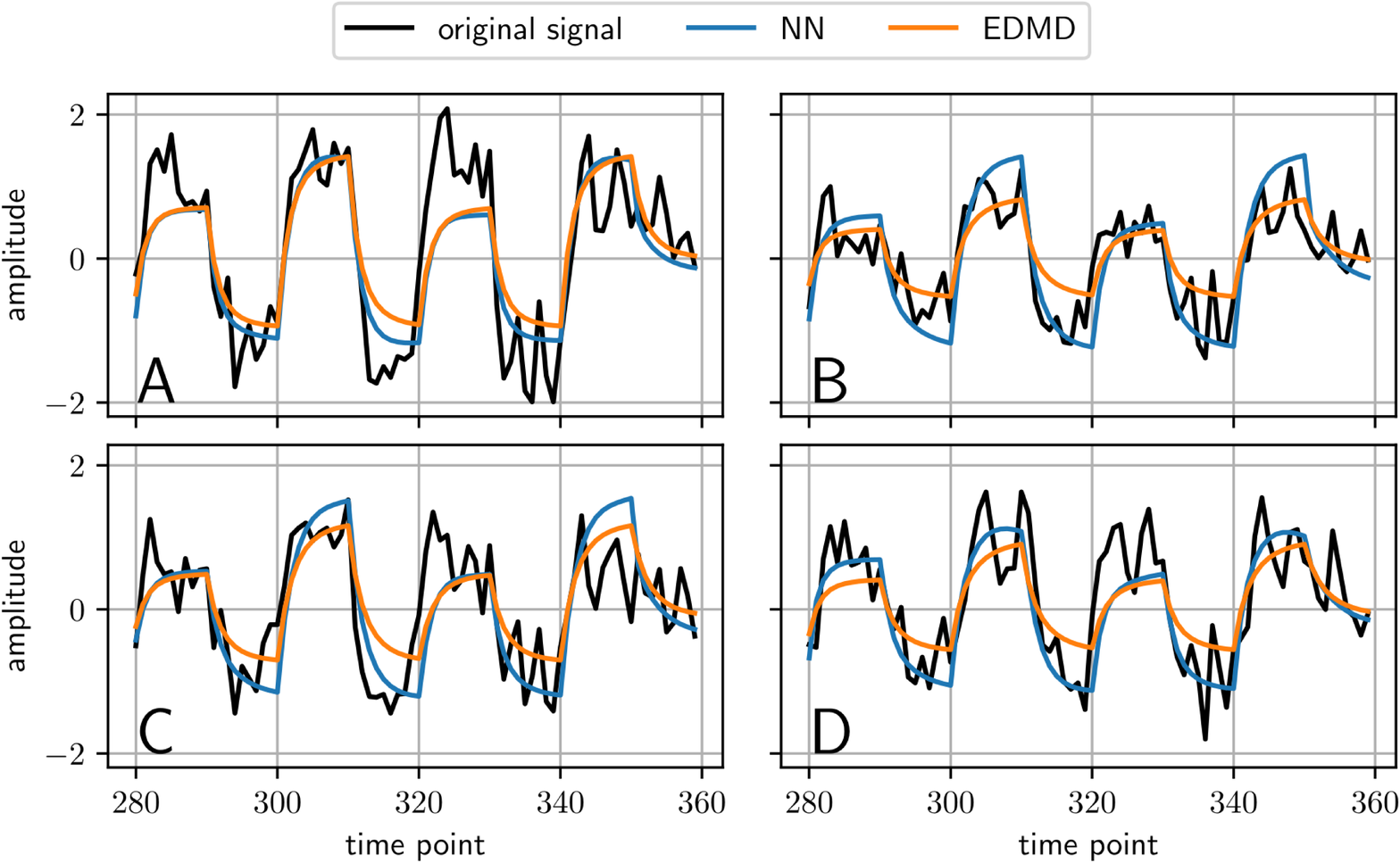}
\caption{Predictions in the amplitude of the BOLD signal in the original fMRI space for four of the most ``activated'' regions as found by the classical GLM methodology  A) left Calcarine Sulcus, B) left Middle Occipital Gyrus, C) right Lingual Gyrus, D) left Superior Occipital Gyrus. The red color marks the actual values of the test set, while the other colors correspond to the predictions based on the proposed methodology.}
\label{fig:figure4}
\end{figure*}

At the final step, the out-of-sample predictions on the manifold made by FNNs and the Koopman operator were ``lifted'' back to the original space using GH and Koopman modes per se, respectively (see \ref{sec:lift}). Specifically, the predicted values were of size $80 \times d=5$ (reduced space) and the lifting of those predictions led to a new set of size $80 \times 111$ (original space). In Figure \ref{fig:figure4}, we depict those predictions overlaid on the unseen (out-of-sample) test data (red color) for the two approaches, namely for the FNN-GH and the Koopman operator. Indicatively, we show the first four regions of interest (ROIs) corresponding to the biggest activated clusters (see Table \ref{tab:table2}), specifically, the left Calcarine sulcus (A), the left Middle Occipital gyrus (B), the right Lingual gyrus (C) and the left Superior Occipital gyrus (D).\par  
In Table \ref{tab:table3}, we present the  errors on the test set for each ROIs as derived for each of the two schemes (FNN- GH and the Koopman operator).\par 
To assess the prediction efficiency of the proposed methodology, we also provide the errors obtained when applying the naive random walk (NRW) model in a direct way (the predictions at time $i+1$ are the last observed values at time $i$). Hence, to predict the values of the DMs at the time point $i+1$ of the test set with the NRW model, we assumed that \textit{we know the actual values of the DMs at the time $i$ of the test set}. The errors are reported in terms of the RMSE and the $L_2$ norm; the smallest errors are marked with bold. In general, the two strategies behave similarly well, outperforming the NRW predictions for all ROIs except one, the Cerebelum Crus 1 on the right hemisphere. This is important, revealing the efficiency of the proposed approach, since as discussed, the predictions of the ROMs were produced iteratively, i.e.  \textit{without any knowledge of the values in the test set}, while the predictions made with the NRW model \textit{were based on the knowledge of the values of the previous points of the test set}.\par 
Here, we should note that predictions made for some brain regions, like the Cuneus of the right hemisphere, might seem to exhibit a poor or spurious predictions. This is due to the fact that in this particular case, the volume of voxels that are ``activated'' (only for the uncorrected threshold of $p<0.001$) is found to contribute only 1.64 \% in the activated cluster (see Table \ref{tab:table2}). Since in this study we try to predict the behaviour of a whole brain region, each time series is averaged over all voxels of the predefined region (here, using AAL). Thus, when trying to predict the ``average'' behaviour of voxels that constitute a region when only a small proportion of them is actually ``activated'', can't be reliable. Therefore, we would normally expect the average signal of such a region to be noisy and without a pattern which could be attributed to the task imposed to the subject (here, the attention to visual motion).

\section{Discussion}
Our work addresses a three-step machine-learning based approach for the modelling of brain activity based on task-depended fMRI data in order to deal with the curse of dimensionality. In the first step, we used parsimonious Diffusion Maps to create a basis that spans a low dimensional subspace, a manifold, in which the information of the brain activity is retained. In a second step, we trained ROMs using this low-dimensional basis, thus relaxing the curse of dimensionality, using both FNNs and the Koopman operator. Finally, we predicted brain activity in the fMRI space by solving the pre-image problem using Geometric Harmonics and the Koopman modes. By doing so, we demonstrated that for the particular benchmark task-experiment, the brain response to the attention to visual motion task is contained in a 5-dimensional manifold. By lifting our predictions to the ambient fMRI space, we found that the key brain regions during attention to visual motion are actually four: the Calcarine Sulcus, the Cuneus, the Superior Occipital Gyrus and Lingual Gyrus. These regions are the same regions reported also in \citep{buchel1998functional} ``re-discovered'' by our methodology.\par 
Our proof-of-concept study is the first to propose such a machine-learning framework to model and importantly accurately predict brain activity from task-based fMRI data, thus facing challenges that are not met when producing controlled experiments from synthetic data produced just by model simulations.\par 
In the numerical experiments detailed above, we show that the performance of the two schemes FNNs-GH and the Koopman operator is comparable, while both outperform the naive random walk whose predictions assume the knowledge of brain activity at the previous time step. Importantly, the comparison between the two schemes reveals that it is not always necessary to construct a non-linear embedding (here, with Diffusion Maps) and then also construct surrogate non-linear models using, e.g., FNNs to predict future states. Effectively, one can utilize the low-frequency truncation of the function space of square-integrable functions ($L^2$) over the original data, as obtained with Diffusion Maps, within the Koopman operator framework, to predict the entire list of coordinate functions in a linear fashion. Conversely, the Neural Network approach treats the embedding coordinates as a base space for the dynamics and predicts them point-wise.\par 
We believe that the proposed methodology may trigger further developments in the field, thus providing the base for a general framework for modelling the dynamics of brain activity, from high-dimensional time series, including also the solution of the source localization problem combining for example simultaneous EEG-fMRI recordings.
\begin{table}
\centering
\caption{ The error evaluated on the test set for each ROI and reconstruction strategy, namely Neural Networks (FNN) and Geometric Harmonics (GH), the Koopman Operator, and the naive Random Walk (NRW) model. The loss is presented in terms of the root mean squared error (RMSE) and the $L_2$ norm (Euclidean distance). For each ROI, the lowest errors are marked with bold.}
\begin{tabular}{ccccccc} 
\toprule
                      & \multicolumn{2}{c}{\textbf{FNN+GH}} & \multicolumn{2}{c}{\textbf{Koopman}} & \multicolumn{2}{c}{\textbf{NRW}}  \\
                      & RMSE  & $L_2$                      & RMSE  & $L_2$                        & RMSE  & $L_2$                        \\
\textbf{Brain Region} &       &                         &       &                           &       &                           \\ 
\midrule
Calcarine L           &\textbf{0.487} & \textbf{4.356}                   & 0.514 & 4.597                     & 0.553 & 4.949                     \\
Calcarine R           & \textbf{0.732} & \textbf{6.545}                   & 0.744 & 6.653                     & 0.772 & 6.902                     \\
Cerebelum 6 L         & \textbf{0.796} & \textbf{7.120}                   & 0.799 & 7.149                     & 0.976 & 8.727                     \\
Cerebelum 6 R         & 0.711 & 6.359                   & \textbf{0.709} & \textbf{6.343}                     & 0.732 & 6.545                     \\
Cerebelum Crus1 L     & 0.690 & 6.176                   & \textbf{0.687} & \textbf{6.142}                     & 0.702 & 6.279                     \\
Cerebelum Crus1 R     & 0.772 & 6.905                   & 0.715 & 6.397                     & \textbf{0.694} & \textbf{6.208}                     \\
Cuneus L              & 0.872 & 7.795                   & \textbf{0.841} & \textbf{7.524}                     & 0.847 & 7.578                     \\
Cuneus R              & 0.965 & 8.628                   & \textbf{0.927} & \textbf{8.291}                     & 1.183 & 10.586                    \\
Fusiform L            & 0.632 & 5.656                   & \textbf{0.601} & \textbf{5.377}                     & 0.658 & 5.886                     \\
Fusiform R            & 0.578 & 5.174                   & \textbf{0.562} & \textbf{5.027}                     & 0.668 & 5.978                     \\
Lingual L             & 0.561 & 5.021                   & \textbf{0.557} & \textbf{4.980}                     & 0.581 & 5.195                     \\
Lingual R             & 0.597 & 5.338                   & \textbf{0.578} & \textbf{5.168}                     & 0.647 & 5.786                     \\
Occipital Inf L       & 0.564 & 5.040                   & 0.\textbf{538} & \textbf{4.810}                     & 0.583 & 5.210                     \\
Occipital Inf R       & 0.902 & 8.064                   & \textbf{0.825} & \textbf{7.379}                     & 1.043 & 9.333                     \\
Occipital Mid L       & \textbf{0.499} & \textbf{4.460}                   & 0.548 & 4.898                     & 0.696 & 6.221                     \\
Occipital Mid R       & 0.755 & 6.754                   & \textbf{0.693} & \textbf{6.200}                     & 0.849 & 7.594                     \\
Occipital Sup L       & \textbf{0.586} & \textbf{5.243}                   & 0.625 & 5.593                     & 0.697 & 6.238                     \\
Occipital Sup R       & \textbf{0.714} & \textbf{6.387}                   & 0.765 & 6.839                     & 0.873 & 7.807                     \\
Parietal Inf R        & \textbf{0.642} & \textbf{5.739}                   & 0.657 & 5.874                     & 0.752 & 6.727                     \\
Postcentral R         & 0.596 & 5.327                   & \textbf{0.590} & \textbf{5.273}                     & 0.711 & 6.359                     \\
\bottomrule
\end{tabular}
\label{tab:table3}
\end{table}

\section{Conflict of Interest}
The authors declare that the research was conducted in the absence of any commercial or financial relationships that could be construed as a potential conflict of interest.

\section{Acknowledgements}
F.D. was funded by the Deutsche Forschungsgemeinschaft (DFG, German Research Foundation) – project no. 468830823.


\begin{thebibliography}{102}
\providecommand{\natexlab}[1]{#1}
\providecommand{\url}[1]{\texttt{#1}}
\expandafter\ifx\csname urlstyle\endcsname\relax
  \providecommand{\doi}[1]{doi: #1}\else
  \providecommand{\doi}{doi: \begingroup \urlstyle{rm}\Url}\fi

\bibitem[Bullmore and Sporns(2009)]{bullmore2009complex}
Ed~Bullmore and Olaf Sporns.
\newblock Complex brain networks: graph theoretical analysis of structural and
  functional systems.
\newblock \emph{Nature reviews neuroscience}, 10\penalty0 (3):\penalty0
  186--198, 2009.

\bibitem[Ermentrout and Terman(2010)]{ermentrout2010mathematical}
Bard Ermentrout and David~H Terman.
\newblock \emph{Mathematical foundations of neuroscience}, volume~35.
\newblock Springer, 2010.

\bibitem[Jorgenson et~al.(2015)Jorgenson, Newsome, Anderson, Bargmann, Brown,
  Deisseroth, Donoghue, Hudson, Ling, MacLeish, et~al.]{jorgenson2015brain}
Lyric~A Jorgenson, William~T Newsome, David~J Anderson, Cornelia~I Bargmann,
  Emery~N Brown, Karl Deisseroth, John~P Donoghue, Kathy~L Hudson, Geoffrey~SF
  Ling, Peter~R MacLeish, et~al.
\newblock The brain initiative: developing technology to catalyse neuroscience
  discovery.
\newblock \emph{Philosophical Transactions of the Royal Society B: Biological
  Sciences}, 370\penalty0 (1668):\penalty0 20140164, 2015.

\bibitem[Siettos and Starke(2016)]{siettos2016multiscale}
Constantinos Siettos and Jens Starke.
\newblock Multiscale modeling of brain dynamics: from single neurons and
  networks to mathematical tools.
\newblock \emph{Wiley Interdisciplinary Reviews: Systems Biology and Medicine},
  8\penalty0 (5):\penalty0 438--458, 2016.

\bibitem[Breakspear(2017)]{breakspear2017dynamic}
Michael Breakspear.
\newblock Dynamic models of large-scale brain activity.
\newblock \emph{Nature neuroscience}, 20\penalty0 (3):\penalty0 340--352, 2017.

\bibitem[Sip et~al.(2023)Sip, Hashemi, Dickscheid, Amunts, Petkoski, and
  Jirsa]{sip2023characterization}
Viktor Sip, Meysam Hashemi, Timo Dickscheid, Katrin Amunts, Spase Petkoski, and
  Viktor Jirsa.
\newblock Characterization of regional differences in resting-state fmri with a
  data-driven network model of brain dynamics.
\newblock \emph{Science Advances}, 9\penalty0 (11):\penalty0 eabq7547, 2023.

\bibitem[Izhikevich(2007)]{izhikevich2007dynamical}
Eugene~M Izhikevich.
\newblock \emph{Dynamical systems in neuroscience}.
\newblock MIT press, 2007.

\bibitem[Deco et~al.(2008)Deco, Jirsa, Robinson, Breakspear, and
  Friston]{deco2008dynamic}
Gustavo Deco, Viktor~K Jirsa, Peter~A Robinson, Michael Breakspear, and Karl
  Friston.
\newblock The dynamic brain: from spiking neurons to neural masses and cortical
  fields.
\newblock \emph{PLoS computational biology}, 4\penalty0 (8):\penalty0 e1000092,
  2008.

\bibitem[Bullmore and Sporns(2012)]{bullmore2012economy}
Ed~Bullmore and Olaf Sporns.
\newblock The economy of brain network organization.
\newblock \emph{Nature reviews neuroscience}, 13\penalty0 (5):\penalty0
  336--349, 2012.

\bibitem[Papo et~al.(2014)Papo, Buld{\'u}, Boccaletti, and
  Bullmore]{papo2014complex}
David Papo, Javier~M Buld{\'u}, Stefano Boccaletti, and Edward~T Bullmore.
\newblock Complex network theory and the brain, 2014.

\bibitem[Hodgkin and Huxley(1952)]{hodgkin1952propagation}
Alan~Lloyd Hodgkin and Andrew~Fielding Huxley.
\newblock Propagation of electrical signals along giant nerve fibres.
\newblock \emph{Proceedings of the Royal Society of London. Series B-Biological
  Sciences}, 140\penalty0 (899):\penalty0 177--183, 1952.

\bibitem[Nelson and Rinzel(1998)]{nelson1998hodgkin}
Mark Nelson and John Rinzel.
\newblock The hodgkin—huxley model.
\newblock In \emph{The book of genesis}, pages 29--49. Springer, 1998.

\bibitem[McCormick et~al.(2007)McCormick, Shu, and Yu]{mccormick2007hodgkin}
David~A McCormick, Yousheng Shu, and Yuguo Yu.
\newblock Hodgkin and huxley model—still standing?
\newblock \emph{Nature}, 445\penalty0 (7123):\penalty0 E1--E2, 2007.

\bibitem[Spiliotis et~al.(2022{\natexlab{a}})Spiliotis, Starke, Franz, Richter,
  and K{\"o}hling]{spiliotis2022deep}
Konstantinos Spiliotis, Jens Starke, Denise Franz, Angelika Richter, and
  R{\"u}diger K{\"o}hling.
\newblock Deep brain stimulation for movement disorder treatment: exploring
  frequency-dependent efficacy in a computational network model.
\newblock \emph{Biological Cybernetics}, 116\penalty0 (1):\penalty0 93--116,
  2022{\natexlab{a}}.

\bibitem[FitzHugh(1955)]{fitzhugh1955mathematical}
Richard FitzHugh.
\newblock Mathematical models of threshold phenomena in the nerve membrane.
\newblock \emph{The bulletin of mathematical biophysics}, 17:\penalty0
  257--278, 1955.

\bibitem[Nagumo et~al.(1962)Nagumo, Arimoto, and Yoshizawa]{nagumo1962active}
Jinichi Nagumo, Suguru Arimoto, and Shuji Yoshizawa.
\newblock An active pulse transmission line simulating nerve axon.
\newblock \emph{Proceedings of the IRE}, 50\penalty0 (10):\penalty0 2061--2070,
  1962.

\bibitem[Kopell and Ermentrout(1986)]{kopell1986symmetry}
Nancy Kopell and GB~Ermentrout.
\newblock Symmetry and phaselocking in chains of weakly coupled oscillators.
\newblock \emph{Communications on Pure and Applied Mathematics}, 39\penalty0
  (5):\penalty0 623--660, 1986.

\bibitem[Laing(2017)]{laing2017phase}
Carlo~R Laing.
\newblock Phase oscillator network models of brain dynamics.
\newblock \emph{Computational models of brain and behavior}, pages 505--517,
  2017.

\bibitem[Skardal and Arenas(2020)]{skardal2020higher}
Per~Sebastian Skardal and Alex Arenas.
\newblock Higher order interactions in complex networks of phase oscillators
  promote abrupt synchronization switching.
\newblock \emph{Communications Physics}, 3\penalty0 (1):\penalty0 218, 2020.

\bibitem[Segneri et~al.(2020)Segneri, Bi, Olmi, and Torcini]{segneri2020theta}
Marco Segneri, Hongjie Bi, Simona Olmi, and Alessandro Torcini.
\newblock Theta-nested gamma oscillations in next generation neural mass
  models.
\newblock \emph{Frontiers in computational neuroscience}, 14:\penalty0 47,
  2020.

\bibitem[Taher et~al.(2020)Taher, Torcini, and Olmi]{taher2020exact}
Halgurd Taher, Alessandro Torcini, and Simona Olmi.
\newblock Exact neural mass model for synaptic-based working memory.
\newblock \emph{PLOS Computational Biology}, 16\penalty0 (12):\penalty0
  e1008533, 2020.

\bibitem[Buxton et~al.(1998)Buxton, Wong, and Frank]{buxton1998dynamics}
Richard~B Buxton, Eric~C Wong, and Lawrence~R Frank.
\newblock Dynamics of blood flow and oxygenation changes during brain
  activation: the balloon model.
\newblock \emph{Magnetic resonance in medicine}, 39\penalty0 (6):\penalty0
  855--864, 1998.

\bibitem[Penny et~al.(2004{\natexlab{a}})Penny, Stephan, Mechelli, and
  Friston]{penny2004modelling}
Will~D Penny, Klaas~E Stephan, Andrea Mechelli, and Karl~J Friston.
\newblock Modelling functional integration: a comparison of structural equation
  and dynamic causal models.
\newblock \emph{Neuroimage}, 23:\penalty0 S264--S274, 2004{\natexlab{a}}.

\bibitem[Heitmann et~al.(2018)Heitmann, Aburn, and
  Breakspear]{heitmann2018brain}
Stewart Heitmann, Matthew~J Aburn, and Michael Breakspear.
\newblock The brain dynamics toolbox for matlab.
\newblock \emph{Neurocomputing}, 315:\penalty0 82--88, 2018.

\bibitem[Beckmann and Smith(2004)]{beckmann2004probabilistic}
Christian~F Beckmann and Stephen~M Smith.
\newblock Probabilistic independent component analysis for functional magnetic
  resonance imaging.
\newblock \emph{IEEE transactions on medical imaging}, 23\penalty0
  (2):\penalty0 137--152, 2004.

\bibitem[Seth(2010)]{seth2010matlab}
Anil~K Seth.
\newblock A matlab toolbox for granger causal connectivity analysis.
\newblock \emph{Journal of neuroscience methods}, 186\penalty0 (2):\penalty0
  262--273, 2010.

\bibitem[Seth et~al.(2015)Seth, Barrett, and Barnett]{seth2015granger}
Anil~K Seth, Adam~B Barrett, and Lionel Barnett.
\newblock Granger causality analysis in neuroscience and neuroimaging.
\newblock \emph{Journal of Neuroscience}, 35\penalty0 (8):\penalty0 3293--3297,
  2015.

\bibitem[Friston et~al.(2013)Friston, Moran, and Seth]{friston2013analysing}
Karl Friston, Rosalyn Moran, and Anil~K Seth.
\newblock Analysing connectivity with granger causality and dynamic causal
  modelling.
\newblock \emph{Current opinion in neurobiology}, 23\penalty0 (2):\penalty0
  172--178, 2013.

\bibitem[Protopapa et~al.(2014)Protopapa, Siettos, Evdokimidis, and
  Smyrnis]{protopapa2014granger}
Foteini Protopapa, Constantinos~I Siettos, Ioannis Evdokimidis, and Nikolaos
  Smyrnis.
\newblock Granger causality analysis reveals distinct spatio-temporal
  connectivity patterns in motor and perceptual visuo-spatial working memory.
\newblock \emph{Frontiers in computational neuroscience}, 8:\penalty0 146,
  2014.

\bibitem[Kugiumtzis and Kimiskidis(2015)]{kugiumtzis2015direct}
Dimitris Kugiumtzis and Vasilios~K Kimiskidis.
\newblock Direct causal networks for the study of transcranial magnetic
  stimulation effects on focal epileptiform discharges.
\newblock \emph{International Journal of Neural Systems}, 25\penalty0
  (05):\penalty0 1550006, 2015.

\bibitem[Protopapa et~al.(2016)Protopapa, Siettos, Myatchin, and
  Lagae]{protopapa2016children}
Foteini Protopapa, Constantinos~I Siettos, Ivan Myatchin, and Lieven Lagae.
\newblock Children with well controlled epilepsy possess different
  spatio-temporal patterns of causal network connectivity during a visual
  working memory task.
\newblock \emph{Cognitive neurodynamics}, 10:\penalty0 99--111, 2016.

\bibitem[Kugiumtzis et~al.(2017)Kugiumtzis, Koutlis, Tsimpiris, and
  Kimiskidis]{kugiumtzis2017dynamics}
Dimitris Kugiumtzis, Christos Koutlis, Alkiviadis Tsimpiris, and Vasilios~K
  Kimiskidis.
\newblock Dynamics of epileptiform discharges induced by transcranial magnetic
  stimulation in genetic generalized epilepsy.
\newblock \emph{International Journal of Neural Systems}, 27\penalty0
  (07):\penalty0 1750037, 2017.

\bibitem[Almpanis and Siettos(2020)]{almpanis2020construction}
Evangelos Almpanis and Constantinos Siettos.
\newblock Construction of functional brain connectivity networks from fmri data
  with driving and modulatory inputs: an extended conditional granger causality
  approach.
\newblock \emph{AIMS neuroscience}, 7\penalty0 (2):\penalty0 66, 2020.

\bibitem[Mormann et~al.(2000)Mormann, Lehnertz, David, and
  Elger]{mormann2000mean}
Florian Mormann, Klaus Lehnertz, Peter David, and Christian~E Elger.
\newblock Mean phase coherence as a measure for phase synchronization and its
  application to the eeg of epilepsy patients.
\newblock \emph{Physica D: Nonlinear Phenomena}, 144\penalty0 (3-4):\penalty0
  358--369, 2000.

\bibitem[Rudrauf et~al.(2006)Rudrauf, Douiri, Kovach, Lachaux, Cosmelli,
  Chavez, Adam, Renault, Martinerie, and Le~Van~Quyen]{rudrauf2006frequency}
David Rudrauf, Abdel Douiri, Christopher Kovach, Jean-Philippe Lachaux, Diego
  Cosmelli, Mario Chavez, Claude Adam, Bernard Renault, Jacques Martinerie, and
  Michel Le~Van~Quyen.
\newblock Frequency flows and the time-frequency dynamics of multivariate phase
  synchronization in brain signals.
\newblock \emph{Neuroimage}, 31\penalty0 (1):\penalty0 209--227, 2006.

\bibitem[Jirsa and M{\"u}ller(2013)]{jirsa2013cross}
Viktor Jirsa and Viktor M{\"u}ller.
\newblock Cross-frequency coupling in real and virtual brain networks.
\newblock \emph{Frontiers in computational neuroscience}, 7:\penalty0 78, 2013.

\bibitem[Zakharova et~al.(2014)Zakharova, Kapeller, and
  Sch{\"o}ll]{zakharova2014chimera}
Anna Zakharova, Marie Kapeller, and Eckehard Sch{\"o}ll.
\newblock Chimera death: Symmetry breaking in dynamical networks.
\newblock \emph{Physical review letters}, 112\penalty0 (15):\penalty0 154101,
  2014.

\bibitem[Mylonas et~al.(2016)Mylonas, Siettos, Evdokimidis, Papanicolaou, and
  Smyrnis]{mylonas2016modular}
DS~Mylonas, CI~Siettos, I~Evdokimidis, AC~Papanicolaou, and N~Smyrnis.
\newblock Modular patterns of phase desynchronization networks during a simple
  visuomotor task.
\newblock \emph{Brain topography}, 29:\penalty0 118--129, 2016.

\bibitem[Sch{\"o}ll(2022)]{scholl2022partial}
Eckehard Sch{\"o}ll.
\newblock Partial synchronization patterns in brain networks.
\newblock \emph{Europhysics Letters}, 136\penalty0 (1):\penalty0 18001, 2022.

\bibitem[Kopell and Ermentrout(2004)]{kopell2004chemical}
Nancy Kopell and Bard Ermentrout.
\newblock Chemical and electrical synapses perform complementary roles in the
  synchronization of interneuronal networks.
\newblock \emph{Proceedings of the National Academy of Sciences}, 101\penalty0
  (43):\penalty0 15482--15487, 2004.

\bibitem[Spiliotis et~al.(2022{\natexlab{b}})Spiliotis, Butenko, van Rienen,
  Starke, and K{\"o}hling]{spiliotis2022complex}
Konstantinos Spiliotis, Konstantin Butenko, Ursula van Rienen, Jens Starke, and
  R{\"u}diger K{\"o}hling.
\newblock Complex network measures reveal optimal targets for deep brain
  stimulation and identify clusters of collective brain dynamics.
\newblock \emph{Frontiers in Physics}, page 1032, 2022{\natexlab{b}}.

\bibitem[Petkoski and Jirsa(2019)]{petkoski2019transmission}
Spase Petkoski and Viktor~K Jirsa.
\newblock Transmission time delays organize the brain network synchronization.
\newblock \emph{Philosophical Transactions of the Royal Society A},
  377\penalty0 (2153):\penalty0 20180132, 2019.

\bibitem[Grossberg and Merrill(1992)]{grossberg1992neural}
Stephen Grossberg and John~WL Merrill.
\newblock A neural network model of adaptively timed reinforcement learning and
  hippocampal dynamics.
\newblock \emph{Cognitive brain research}, 1\penalty0 (1):\penalty0 3--38,
  1992.

\bibitem[Niv(2009)]{niv2009reinforcement}
Yael Niv.
\newblock Reinforcement learning in the brain.
\newblock \emph{Journal of Mathematical Psychology}, 53\penalty0 (3):\penalty0
  139--154, 2009.

\bibitem[Richiardi et~al.(2013)Richiardi, Achard, Bunke, and Van
  De~Ville]{richiardi2013machine}
Jonas Richiardi, Sophie Achard, Horst Bunke, and Dimitri Van De~Ville.
\newblock Machine learning with brain graphs: predictive modeling approaches
  for functional imaging in systems neuroscience.
\newblock \emph{IEEE Signal processing magazine}, 30\penalty0 (3):\penalty0
  58--70, 2013.

\bibitem[Suk et~al.(2016)Suk, Wee, Lee, and Shen]{suk2016state}
Heung-Il Suk, Chong-Yaw Wee, Seong-Whan Lee, and Dinggang Shen.
\newblock State-space model with deep learning for functional dynamics
  estimation in resting-state fmri.
\newblock \emph{NeuroImage}, 129:\penalty0 292--307, 2016.

\bibitem[Gholami~Doborjeh et~al.(2018)Gholami~Doborjeh, Kasabov,
  Gholami~Doborjeh, and Sumich]{gholami2018modelling}
Zohreh Gholami~Doborjeh, Nikola Kasabov, Maryam Gholami~Doborjeh, and Alexander
  Sumich.
\newblock Modelling peri-perceptual brain processes in a deep learning spiking
  neural network architecture.
\newblock \emph{Scientific reports}, 8\penalty0 (1):\penalty0 8912, 2018.

\bibitem[Sun et~al.(2022)Sun, Sohrabpour, Worrell, and He]{sun2022deep}
Rui Sun, Abbas Sohrabpour, Gregory~A Worrell, and Bin He.
\newblock Deep neural networks constrained by neural mass models improve
  electrophysiological source imaging of spatiotemporal brain dynamics.
\newblock \emph{Proceedings of the National Academy of Sciences}, 119\penalty0
  (31):\penalty0 e2201128119, 2022.

\bibitem[Phinyomark et~al.(2017)Phinyomark, Ibanez-Marcelo, and
  Petri]{phinyomark2017resting}
Angkoon Phinyomark, Esther Ibanez-Marcelo, and Giovanni Petri.
\newblock Resting-state fmri functional connectivity: Big data preprocessing
  pipelines and topological data analysis.
\newblock \emph{IEEE Transactions on Big Data}, 3\penalty0 (4):\penalty0
  415--428, 2017.

\bibitem[Altman and Krzywinski(2018)]{altman2018curse}
Naomi Altman and Martin Krzywinski.
\newblock The curse (s) of dimensionality.
\newblock \emph{Nat Methods}, 15\penalty0 (6):\penalty0 399--400, 2018.

\bibitem[Penny et~al.(2011)Penny, Friston, Ashburner, Kiebel, and
  Nichols]{penny2011statistical}
William~D Penny, Karl~J Friston, John~T Ashburner, Stefan~J Kiebel, and
  Thomas~E Nichols.
\newblock \emph{Statistical parametric mapping: the analysis of functional
  brain images}.
\newblock Elsevier, 2011.

\bibitem[Madhyastha et~al.(2018)Madhyastha, Peverill, Koh, McCabe, Flournoy,
  Mills, King, Pfeifer, and McLaughlin]{madhyastha2018current}
Tara Madhyastha, Matthew Peverill, Natalie Koh, Connor McCabe, John Flournoy,
  Kate Mills, Kevin King, Jennifer Pfeifer, and Katie~A McLaughlin.
\newblock Current methods and limitations for longitudinal fmri analysis across
  development.
\newblock \emph{Developmental cognitive neuroscience}, 33:\penalty0 118--128,
  2018.

\bibitem[Belkin and Niyogi(2003)]{belkin2003laplacian}
Mikhail Belkin and Partha Niyogi.
\newblock Laplacian eigenmaps for dimensionality reduction and data
  representation.
\newblock \emph{Neural computation}, 15\penalty0 (6):\penalty0 1373--1396,
  2003.

\bibitem[Coifman and Lafon(2006{\natexlab{a}})]{coifman2006diffusion}
Ronald~R Coifman and St{\'e}phane Lafon.
\newblock Diffusion maps.
\newblock \emph{Applied and computational harmonic analysis}, 21\penalty0
  (1):\penalty0 5--30, 2006{\natexlab{a}}.

\bibitem[Ansuini et~al.(2019)Ansuini, Laio, Macke, and
  Zoccolan]{ansuini2019intrinsic}
Alessio Ansuini, Alessandro Laio, Jakob~H Macke, and Davide Zoccolan.
\newblock Intrinsic dimension of data representations in deep neural networks.
\newblock \emph{Advances in Neural Information Processing Systems}, 32, 2019.

\bibitem[Gallos et~al.(2021{\natexlab{a}})Gallos, Galaris, and
  Siettos]{gallos2021construction}
Ioannis~K Gallos, Evangelos Galaris, and Constantinos~I Siettos.
\newblock Construction of embedded fmri resting-state functional connectivity
  networks using manifold learning.
\newblock \emph{Cognitive neurodynamics}, 15\penalty0 (4):\penalty0 585--608,
  2021{\natexlab{a}}.

\bibitem[Qiu et~al.(2015)Qiu, Lee, Tan, and Chung]{qiu2015manifold}
Anqi Qiu, Annie Lee, Mingzhen Tan, and Moo~K Chung.
\newblock Manifold learning on brain functional networks in aging.
\newblock \emph{Medical image analysis}, 20\penalty0 (1):\penalty0 52--60,
  2015.

\bibitem[Roweis and Saul(2000)]{roweis2000nonlinear}
Sam~T Roweis and Lawrence~K Saul.
\newblock Nonlinear dimensionality reduction by locally linear embedding.
\newblock \emph{science}, 290\penalty0 (5500):\penalty0 2323--2326, 2000.

\bibitem[Pospelov et~al.(2021)Pospelov, Tetereva, Martynova, and
  Anokhin]{pospelov2021laplacian}
Nikita Pospelov, Alina Tetereva, Olga Martynova, and Konstantin Anokhin.
\newblock The laplacian eigenmaps dimensionality reduction of fmri data for
  discovering stimulus-induced changes in the resting-state brain activity.
\newblock \emph{Neuroimage: Reports}, 1\penalty0 (3):\penalty0 100035, 2021.

\bibitem[Gao et~al.(2021)Gao, Mishne, and Scheinost]{gao2021nonlinear}
Siyuan Gao, Gal Mishne, and Dustin Scheinost.
\newblock Nonlinear manifold learning in functional magnetic resonance imaging
  uncovers a low-dimensional space of brain dynamics.
\newblock \emph{Human brain mapping}, 42\penalty0 (14):\penalty0 4510--4524,
  2021.

\bibitem[Gallos et~al.(2021{\natexlab{b}})Gallos, Gkiatis, Matsopoulos, and
  Siettos]{gallos2021isomap}
Ioannis~K Gallos, Kostakis Gkiatis, George~K Matsopoulos, and Constantinos
  Siettos.
\newblock Isomap and machine learning algorithms for the construction of
  embedded functional connectivity networks of anatomically separated brain
  regions from resting state fmri data of patients with schizophrenia.
\newblock \emph{AIMS neuroscience}, 8\penalty0 (2):\penalty0 295,
  2021{\natexlab{b}}.

\bibitem[Gallos et~al.(2021{\natexlab{c}})Gallos, Mantonakis, Spilioti,
  Kattoulas, Savvidou, Anyfandi, Karavasilis, Kelekis, Smyrnis, and
  Siettos]{gallos2021relation}
IK~Gallos, L~Mantonakis, E~Spilioti, E~Kattoulas, E~Savvidou, E~Anyfandi,
  E~Karavasilis, N~Kelekis, N~Smyrnis, and CI~Siettos.
\newblock The relation of integrated psychological therapy to resting state
  functional brain connectivity networks in patients with schizophrenia.
\newblock \emph{Psychiatry Research}, page 114270, 2021{\natexlab{c}}.

\bibitem[Cox and Cox(2008)]{cox2008multidimensional}
Michael~AA Cox and Trevor~F Cox.
\newblock Multidimensional scaling.
\newblock In \emph{Handbook of data visualization}, pages 315--347. Springer,
  2008.

\bibitem[Tenenbaum et~al.(2000)Tenenbaum, De~Silva, and
  Langford]{tenenbaum2000global}
Joshua~B Tenenbaum, Vin De~Silva, and John~C Langford.
\newblock A global geometric framework for nonlinear dimensionality reduction.
\newblock \emph{science}, 290\penalty0 (5500):\penalty0 2319--2323, 2000.

\bibitem[Nadler et~al.(2006)Nadler, Lafon, Kevrekidis, and
  Coifman]{nadler2006diffusion}
Boaz Nadler, Stephane Lafon, Ioannis Kevrekidis, and Ronald~R Coifman.
\newblock Diffusion maps, spectral clustering and eigenfunctions of
  fokker-planck operators.
\newblock In \emph{Advances in neural information processing systems}, pages
  955--962, 2006.

\bibitem[Nadler et~al.(2008)Nadler, Lafon, Coifman, and
  Kevrekidis]{nadler2008diffusion}
Boaz Nadler, Stephane Lafon, Ronald Coifman, and Ioannis~G. Kevrekidis.
\newblock Diffusion maps - a probabilistic interpretation for spectral
  embedding and clustering algorithms.
\newblock In Alexander~N. Gorban, Bal{\'a}zs K{\'e}gl, Donald~C. Wunsch, and
  Andrei~Y. Zinovyev, editors, \emph{Principal Manifolds for Data Visualization
  and Dimension Reduction}, pages 238--260, Berlin, Heidelberg, 2008. Springer
  Berlin Heidelberg.
\newblock ISBN 978-3-540-73750-6.

\bibitem[Thiem et~al.(2020)Thiem, Kooshkbaghi, Bertalan, Laing, and
  Kevrekidis]{thiem2020emergent}
Thomas~N Thiem, Mahdi Kooshkbaghi, Tom Bertalan, Carlo~R Laing, and Ioannis~G
  Kevrekidis.
\newblock Emergent spaces for coupled oscillators.
\newblock \emph{Frontiers in computational neuroscience}, 14:\penalty0 36,
  2020.

\bibitem[Mezi{\'c}(2013)]{mezic2013analysis}
Igor Mezi{\'c}.
\newblock Analysis of fluid flows via spectral properties of the koopman
  operator.
\newblock \emph{Annual Review of Fluid Mechanics}, 45:\penalty0 357--378, 2013.

\bibitem[Williams et~al.(2015{\natexlab{a}})Williams, Kevrekidis, and
  Rowley]{williams2015data}
Matthew~O Williams, Ioannis~G Kevrekidis, and Clarence~W Rowley.
\newblock A data--driven approximation of the koopman operator: Extending
  dynamic mode decomposition.
\newblock \emph{Journal of Nonlinear Science}, 25:\penalty0 1307--1346,
  2015{\natexlab{a}}.

\bibitem[Brunton et~al.(2016)Brunton, Brunton, Proctor, and
  Kutz]{brunton2016koopman}
Steven~L Brunton, Bingni~W Brunton, Joshua~L Proctor, and J~Nathan Kutz.
\newblock Koopman invariant subspaces and finite linear representations of
  nonlinear dynamical systems for control.
\newblock \emph{PloS one}, 11\penalty0 (2):\penalty0 e0150171, 2016.

\bibitem[Li et~al.(2017)Li, Dietrich, Bollt, and Kevrekidis]{li2017extended}
Qianxiao Li, Felix Dietrich, Erik~M Bollt, and Ioannis~G Kevrekidis.
\newblock Extended dynamic mode decomposition with dictionary learning: A
  data-driven adaptive spectral decomposition of the koopman operator.
\newblock \emph{Chaos: An Interdisciplinary Journal of Nonlinear Science},
  27\penalty0 (10):\penalty0 103111, 2017.

\bibitem[Bollt et~al.(2018)Bollt, Li, Dietrich, and
  Kevrekidis]{bollt2018matching}
Erik~M Bollt, Qianxiao Li, Felix Dietrich, and Ioannis Kevrekidis.
\newblock On matching, and even rectifying, dynamical systems through koopman
  operator eigenfunctions.
\newblock \emph{SIAM Journal on Applied Dynamical Systems}, 17\penalty0
  (2):\penalty0 1925--1960, 2018.

\bibitem[Dietrich et~al.(2020{\natexlab{a}})Dietrich, Thiem, and
  Kevrekidis]{dietrich2020koopman}
Felix Dietrich, Thomas~N Thiem, and Ioannis~G Kevrekidis.
\newblock On the koopman operator of algorithms.
\newblock \emph{SIAM Journal on Applied Dynamical Systems}, 19\penalty0
  (2):\penalty0 860--885, 2020{\natexlab{a}}.

\bibitem[Lehmberg et~al.(2021)Lehmberg, Dietrich, and
  K{\"o}ster]{lehmberg2021modeling}
Daniel Lehmberg, Felix Dietrich, and Gerta K{\"o}ster.
\newblock Modeling melburnians—using the koopman operator to gain insight
  into crowd dynamics.
\newblock \emph{Transportation Research Part C: Emerging Technologies},
  133:\penalty0 103437, 2021.

\bibitem[Coifman and Lafon(2006{\natexlab{b}})]{coifman2006geometric}
Ronald~R Coifman and St{\'e}phane Lafon.
\newblock Geometric harmonics: a novel tool for multiscale out-of-sample
  extension of empirical functions.
\newblock \emph{Applied and Computational Harmonic Analysis}, 21\penalty0
  (1):\penalty0 31--52, 2006{\natexlab{b}}.

\bibitem[Dsilva et~al.(2013)Dsilva, Talmon, Rabin, Coifman, and
  Kevrekidis]{dsilva2013nonlinear}
Carmeline~J Dsilva, Ronen Talmon, Neta Rabin, Ronald~R Coifman, and Ioannis~G
  Kevrekidis.
\newblock Nonlinear intrinsic variables and state reconstruction in multiscale
  simulations.
\newblock \emph{The Journal of chemical physics}, 139\penalty0 (18):\penalty0
  11B608\_1, 2013.

\bibitem[Papaioannou et~al.(2021)Papaioannou, Talmon, di~Serafino, Kevrekidis,
  and Siettos]{papaioannou2021time}
Panagiotis Papaioannou, Ronen Talmon, Daniela di~Serafino, Ioannis Kevrekidis,
  and Constantinos Siettos.
\newblock Time series forecasting using manifold learning.
\newblock \emph{arXiv preprint arXiv:2110.03625}, 2021.

\bibitem[Evangelou et~al.(2023)Evangelou, Dietrich, Chiavazzo, Lehmberg, Meila,
  and Kevrekidis]{evangelou2022double}
Nikolaos Evangelou, Felix Dietrich, Eliodoro Chiavazzo, Daniel Lehmberg, Marina
  Meila, and Ioannis~G. Kevrekidis.
\newblock Double diffusion maps and their latent harmonics for scientific
  computations in latent space.
\newblock \emph{Journal of Computational Physics}, page 112072, 2023.

\bibitem[B{\"u}chel and Friston(1997)]{buchel1997modulation}
Christian B{\"u}chel and Karl~J Friston.
\newblock Modulation of connectivity in visual pathways by attention: cortical
  interactions evaluated with structural equation modelling and fmri.
\newblock \emph{Cerebral Cortex (New York, NY: 1991)}, 7\penalty0 (8):\penalty0
  768--778, 1997.

\bibitem[Friston et~al.(2003)Friston, Harrison, and Penny]{friston2003dynamic}
Karl~J Friston, Lee Harrison, and Will Penny.
\newblock Dynamic causal modelling.
\newblock \emph{Neuroimage}, 19\penalty0 (4):\penalty0 1273--1302, 2003.

\bibitem[Penny et~al.(2004{\natexlab{b}})Penny, Stephan, Mechelli, and
  Friston]{penny2004comparing}
Will~D Penny, Klaas~E Stephan, Andrea Mechelli, and Karl~J Friston.
\newblock Comparing dynamic causal models.
\newblock \emph{Neuroimage}, 22\penalty0 (3):\penalty0 1157--1172,
  2004{\natexlab{b}}.

\bibitem[Tzourio-Mazoyer et~al.(2002)Tzourio-Mazoyer, Landeau, Papathanassiou,
  Crivello, Etard, Delcroix, Mazoyer, and Joliot]{tzourio2002automated}
Nathalie Tzourio-Mazoyer, Brigitte Landeau, Dimitri Papathanassiou, Fabrice
  Crivello, Olivier Etard, Nicolas Delcroix, Bernard Mazoyer, and Marc Joliot.
\newblock Automated anatomical labeling of activations in spm using a
  macroscopic anatomical parcellation of the mni mri single-subject brain.
\newblock \emph{Neuroimage}, 15\penalty0 (1):\penalty0 273--289, 2002.

\bibitem[Williams et~al.(2015{\natexlab{b}})Williams, Kevrekidis, and
  Rowley]{williams-2015b}
Matthew~O. Williams, Ioannis~G. Kevrekidis, and Clarence~W. Rowley.
\newblock A {{Data-Driven Approximation}} of the {{Koopman Operator}}:
  {{Extending Dynamic Mode Decomposition}}.
\newblock \emph{Journal of Nonlinear Science}, 25\penalty0 (6):\penalty0
  1307--1346, December 2015{\natexlab{b}}.
\newblock ISSN 0938-8974, 1432-1467.
\newblock \doi{10.1007/s00332-015-9258-5}.

\bibitem[Friston et~al.(1994)Friston, Holmes, Worsley, Poline, Frith, and
  Frackowiak]{friston1994statistical}
Karl~J Friston, Andrew~P Holmes, Keith~J Worsley, J-P Poline, Chris~D Frith,
  and Richard~SJ Frackowiak.
\newblock Statistical parametric maps in functional imaging: a general linear
  approach.
\newblock \emph{Human brain mapping}, 2\penalty0 (4):\penalty0 189--210, 1994.

\bibitem[Friston et~al.(1995)Friston, Holmes, Poline, Grasby, Williams,
  Frackowiak, and Turner]{friston1995analysis}
Karl~J Friston, Andrew~P Holmes, JB~Poline, PJ~Grasby, SCR Williams, Richard~SJ
  Frackowiak, and Robert Turner.
\newblock Analysis of fmri time-series revisited.
\newblock \emph{Neuroimage}, 2\penalty0 (1):\penalty0 45--53, 1995.

\bibitem[Worsley and Friston(1995)]{worsley1995analysis}
Keith~J Worsley and Karl~J Friston.
\newblock Analysis of fmri time-series revisited—again.
\newblock \emph{Neuroimage}, 2\penalty0 (3):\penalty0 173--181, 1995.

\bibitem[B{\"u}chel et~al.(1998)B{\"u}chel, Josephs, Rees, Turner, Frith, and
  Friston]{buchel1998functional}
Christian B{\"u}chel, Oliver Josephs, Geraint Rees, Robert Turner, Chris~D
  Frith, and Karl~J Friston.
\newblock The functional anatomy of attention to visual motion. a functional
  mri study.
\newblock \emph{Brain: a journal of neurology}, 121\penalty0 (7):\penalty0
  1281--1294, 1998.

\bibitem[Dsilva et~al.(2018)Dsilva, Talmon, Coifman, and
  Kevrekidis]{dsilva2018parsimonious}
Carmeline~J Dsilva, Ronen Talmon, Ronald~R Coifman, and Ioannis~G Kevrekidis.
\newblock Parsimonious representation of nonlinear dynamical systems through
  manifold learning: A chemotaxis case study.
\newblock \emph{Applied and Computational Harmonic Analysis}, 44\penalty0
  (3):\penalty0 759--773, 2018.

\bibitem[Lee et~al.(2020)Lee, Kooshkbaghi, Spiliotis, Siettos, and
  Kevrekidis]{lee2020coarse}
Seungjoon Lee, Mahdi Kooshkbaghi, Konstantinos Spiliotis, Constantinos~I
  Siettos, and Ioannis~G Kevrekidis.
\newblock Coarse-scale pdes from fine-scale observations via machine learning.
\newblock \emph{Chaos: An Interdisciplinary Journal of Nonlinear Science},
  30\penalty0 (1):\penalty0 013141, 2020.

\bibitem[Galaris et~al.(2022)Galaris, Fabiani, Gallos, Kevrekidis, and
  Siettos]{galaris2022numerical}
Evangelos Galaris, Gianluca Fabiani, Ioannis Gallos, Ioannis Kevrekidis, and
  Constantinos Siettos.
\newblock Numerical bifurcation analysis of pdes from lattice boltzmann model
  simulations: a parsimonious machine learning approach.
\newblock \emph{Journal of Scientific Computing}, 92\penalty0 (2):\penalty0 34,
  2022.

\bibitem[Krogh and Hertz(1992)]{krogh1992simple}
Anders Krogh and John~A Hertz.
\newblock A simple weight decay can improve generalization.
\newblock In \emph{Advances in neural information processing systems}, pages
  950--957, 1992.

\bibitem[Mezi{\'c}(2005)]{mezic-2005}
Igor Mezi{\'c}.
\newblock Spectral {{Properties}} of {{Dynamical Systems}}, {{Model Reduction}}
  and {{Decompositions}}.
\newblock \emph{Nonlinear Dynamics}, 41\penalty0 (1):\penalty0 309--325, August
  2005.
\newblock ISSN 1573-269X.
\newblock \doi{10.1007/s11071-005-2824-x}.

\bibitem[Budi{\v s}i{\'c} et~al.(2012)Budi{\v s}i{\'c}, Mohr, and
  Mezi{\'c}]{budisic-2012}
Marko Budi{\v s}i{\'c}, Ryan Mohr, and Igor Mezi{\'c}.
\newblock Applied {{Koopmanism}}.
\newblock \emph{Chaos: An Interdisciplinary Journal of Nonlinear Science},
  22:\penalty0 047510, 2012.
\newblock \doi{10.1063/1.4772195}.

\bibitem[Dietrich et~al.(2020{\natexlab{b}})Dietrich, Thiem, and
  Kevrekidis]{dietrich-2020a}
Felix Dietrich, Thomas~N. Thiem, and Ioannis~G. Kevrekidis.
\newblock On the {{Koopman Operator}} of {{Algorithms}}.
\newblock \emph{SIAM Journal on Applied Dynamical Systems}, 19\penalty0
  (2):\penalty0 860--885, January 2020{\natexlab{b}}.
\newblock \doi{10.1137/19m1277059}.

\bibitem[Schmid(2010)]{schmid-2010}
Peter~J. Schmid.
\newblock Dynamic mode decomposition of numerical and experimental data.
\newblock \emph{Journal of Fluid Mechanics}, 656:\penalty0 5--28, 2010.
\newblock \doi{10.1017/s0022112010001217}.

\bibitem[Schmid(2022)]{schmid-2022}
Peter~J. Schmid.
\newblock Dynamic {{Mode Decomposition}} and {{Its Variants}}.
\newblock \emph{Annual Review of Fluid Mechanics}, 54\penalty0 (1):\penalty0
  225--254, January 2022.
\newblock ISSN 0066-4189, 1545-4479.
\newblock \doi{10.1146/annurev-fluid-030121-015835}.

\bibitem[Chiavazzo et~al.(2014)Chiavazzo, Gear, Dsilva, Rabin, and
  Kevrekidis]{chiavazzo2014reduced}
Eliodoro Chiavazzo, Charles~W Gear, Carmeline~J Dsilva, Neta Rabin, and
  Ioannis~G Kevrekidis.
\newblock Reduced models in chemical kinetics via nonlinear data-mining.
\newblock \emph{Processes}, 2\penalty0 (1):\penalty0 112--140, 2014.

\bibitem[Nystr{\"o}m(1929)]{nystrom1929praktische}
Evert~Johannes Nystr{\"o}m.
\newblock \emph{{\"U}ber die praktische aufl{\"o}sung von linearen
  integralgleichungen mit anwendungen auf randwertaufgaben der
  potentialtheorie}.
\newblock Akademische Buchhandlung, 1929.

\bibitem[Patsatzis et~al.(2023)Patsatzis, Russo, Kevrekidis, and
  Siettos]{patsatzis2023data}
Dimitrios~G Patsatzis, Lucia Russo, Ioannis~G Kevrekidis, and Constantinos
  Siettos.
\newblock Data-driven control of agent-based models: An equation/variable-free
  machine learning approach.
\newblock \emph{Journal of Computational Physics}, 478:\penalty0 111953, 2023.

\bibitem[Lehmberg et~al.(2020)Lehmberg, Dietrich, K{\"o}ster, and
  Bungartz]{lehmberg2020datafold}
Daniel Lehmberg, Felix Dietrich, Gerta K{\"o}ster, and Hans-Joachim Bungartz.
\newblock Datafold: data-driven models for point clouds and time series on
  manifolds.
\newblock \emph{Journal of Open Source Software}, 5\penalty0 (51):\penalty0
  2283, 2020.

\bibitem[Pedregosa et~al.(2011)Pedregosa, Varoquaux, Gramfort, Michel, Thirion,
  Grisel, Blondel, Prettenhofer, Weiss, Dubourg, et~al.]{pedregosa2011scikit}
Fabian Pedregosa, Ga{\"e}l Varoquaux, Alexandre Gramfort, Vincent Michel,
  Bertrand Thirion, Olivier Grisel, Mathieu Blondel, Peter Prettenhofer, Ron
  Weiss, Vincent Dubourg, et~al.
\newblock Scikit-learn: Machine learning in python.
\newblock \emph{the Journal of machine Learning research}, 12:\penalty0
  2825--2830, 2011.

\bibitem[Ripley et~al.(2016)Ripley, Venables, and Ripley]{ripley2016package}
Brian Ripley, William Venables, and Maintainer~Brian Ripley.
\newblock Package ‘nnet’.
\newblock \emph{R package version}, 7\penalty0 (3-12):\penalty0 700, 2016.

\end{thebibliography}

\end{document}